\documentclass[12pt,a4paper,twoside,final,notitlepage, leqno]{article}
\usepackage[english]{babel}
\usepackage[T1]{fontenc}  
\usepackage{graphicx}
\usepackage{amsmath,amsthm,epsfig,amsfonts,bbm}   

\def \smb {{\scriptstyle \bullet }}
\newcommand{\monitem}{ \smallskip \noindent $\bullet$ \quad  } 
\newcommand{\moneq}{\vspace*{-6pt} \begin{equation} \displaystyle } 
\newcommand{\moneqstar}{\vspace*{-6pt} \begin{equation*} \displaystyle } 
\newcommand{\monendstar}{\vspace*{-6pt} \end{equation*}   }
\newcommand{\monend}{\vspace*{-6pt} \end{equation}   }

\def\br {\break}

\def\section*#1{}

\def\resume{\if@twocolumn
\section*{R\'esum\'e}
\else \small
\quotation{\bf \it R\'esum\'e \rule[1mm]{1.5mm}{0.2mm}\vspace{0pt}}
\fi}
\def\endresume{\if@twocolumn\else\endquotation\fi}
\def\abstract{\if@twocolumn
\noindent\section*{{\bf Abstract}}
\else \small
\quotation{\noindent \bf {Abstract.} \rule[1mm]{1.5mm}{0.2mm}\vspace{0pt}}
\fi}
\def\endabstract{\if@twocolumn\else\endquotation\fi}

\usepackage{fancyhdr}
\renewcommand{\headrulewidth}{0pt}
\begin{document}

\fancypagestyle{plain}{ \fancyfoot{} \renewcommand{\footrulewidth}{0pt}}
\fancypagestyle{plain}{ \fancyhead{} \renewcommand{\headrulewidth}{0pt}}

\bibliographystyle{alpha}


\title{\bf \LARGE  Some results on energy-conserving \\~ lattice Boltzmann models  \\~ }

\author { { ~  \large    Pierre Lallemand $^1$  and Fran\c{c}ois Dubois $^{2,}$$^3$} \\ ~ \\   
{\it \small $^1$   Centre National de la Recherche Scientifique, Paris, France}\\    
{\it \small $^2$   Conservatoire National des Arts et M\'etiers,  Paris, France}\\  
{\it \small $^3$  Department of Mathematics, University  Paris Sud, Orsay, France.}\\   ~ \\
{ \rm  \small  pierre.lallemand1@free.fr , francois.dubois@math.u-psud.fr . }}

\date  {{  \rm  28  November  2012} 
 \footnote {\rm  \small $\,\,$ Contribution published  in the journal   
 {\it Computers and Mathematics with Applications},  
volume 65, pages 831-844, march 2013,   doi:10.1016/j.camwa.2012.11.007.  
The original communication was 
presented  at the  Eighth International Conference for Mesoscopic   
Methods in Engineering  and Science,  Lyon, 6 July   2011. 
Edition 11 september 2013.  }}

\maketitle
\renewcommand{\baselinestretch}{1.}

  \vspace{-.5cm}  
  
 \bigskip \bigskip  
\noindent {\bf  \large   Abstract  } 
 
\noindent  
We   consider  the  problem of  ``energy conserving''  lattice Boltzmann models. 
A major difficulty observed in previous studies  is the coupling between 
the viscous and thermal waves even at moderate wave numbers. 
We propose a theoretical framework based 
on the  knowledge of the partial equivalent equations of the lattice 
Boltzmann scheme at several orders of precision. 
With the help of linearized models (inviscid and dissipative advective acoustics
and classical acoustics), we suggest natural sets of relations for the parameters   
of lattice Boltzmann schemes. 
The application is proposed for three two-dimensional schemes. 
Numerical test cases for simple linear and nonlinear waves establish
that the main difficulty in the previous contributions can now be
overcome.                                      

 $ $ \\[.4mm]

 $ $ \\[.4mm]

\newpage 
\fancyfoot[C]{\oldstylenums{\thepage}}

~

 \noindent {\bf \large 1) \quad    Introduction  } 
 
 \noindent  
 \noindent  
The derivation of  lattice gas automata taking into account 
the  conservation of mass, momentum and total 
energy has been initially  proposed by McNamara and Alder \cite{MA93}. 
In his contribution that fixed the paradigm of multiple
relaxation times of lattice Boltzmann schemes, 
 d'Humi\`eres \cite{DdH92} presented simulation
of compressible fluids with the presence of strong discontinuities. 
Nevertheless, in order to fit the equilibrium distribution, 
it is necessary to consider lattice Boltzmann models with a large number
of velocities (see {\it e.g.}  Qian \cite{Qi93} and 
 Alexander, Chen and  Sterling \cite{ACS93}).
In the  study of one of us with L.S.~Luo \cite{LL03}, 
it has been established that the  D2Q9 scheme 
(see the  Figure~10 and a detailed description in  Annex~1)
does not allow correctly a variation of sound velocity with the temperature. 
The contribution  \cite{LL03}  enforces the use of higher order stencils
as the  D2Q13 scheme (see   Figure~11 and    Annex~2).

\smallskip   \smallskip  \centerline {                       
\includegraphics[width=.52 \textwidth]  {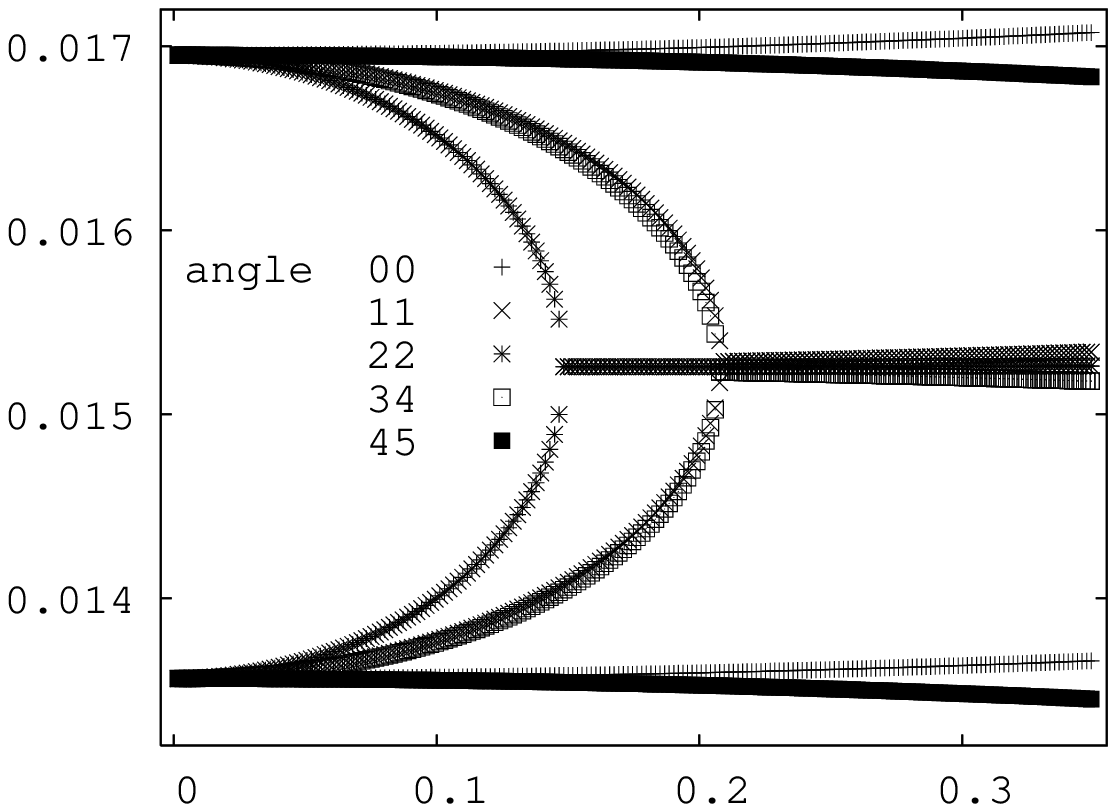} 
\includegraphics[width=.52 \textwidth]  {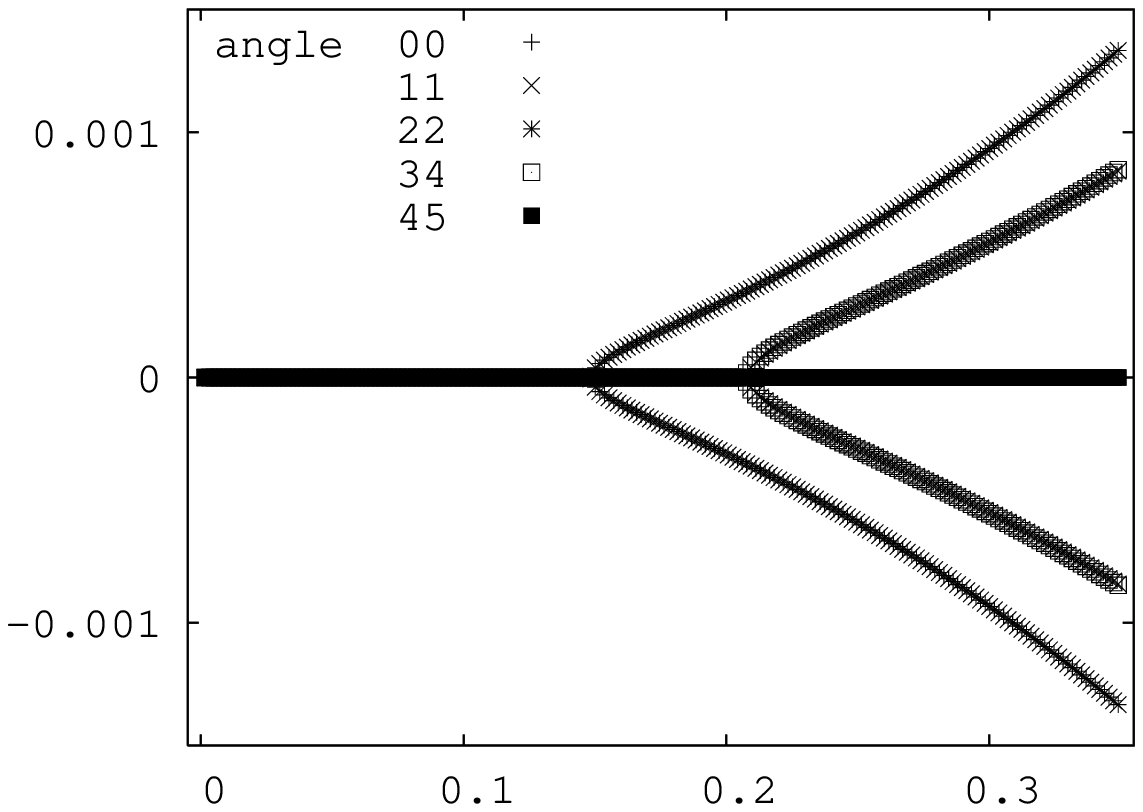} }

\smallskip \noindent  {\bf Figure 1}. \quad  ``Zero point'' experiment 
with the D2Q13 lattice Boltzmann scheme.   The wave vector is in abscissa and the 
eigenvalues of the lattice Boltzmann scheme in y-axis. The real part is on the
left figure and the imaginary part on the right. 
Results for two different angles. 
For a critical wave number, the viscous and thermal modes merge together and the physics
is badly approximated. 
Choice of coefficients defined in the  relations  (\ref{relax}) 
and  (\ref{first-equil-D2Q13})~: 
$ \,  c_1=-1.3 \,,\,  \alpha_2=-25 \,,\, \beta_2=-1.5
 \,,\, \alpha_3=5.5 \,,\,  \beta_3=0, \, $ 
$ s_5 = 1.88 ,\,  s_7 =1.95  ,\,  s_9 =1.60  ,\,  s_{11} =1.75 
,\,  s_{12} =1.05  ,\,   s_{13} =1.35 .   $

\smallskip \smallskip 

\monitem 
A major difficulty observed in  \cite{LL03} is the coupling between 
the viscous and thermal waves at moderate wave numbers. Consider 
a DdQq lattice Boltzmann scheme with discrete velocities $ \, \xi_j \, $
(see (\ref{vit-schema}),  (\ref{vit-d2q9}),    (\ref{vit-d2q13}) and  (\ref{vit-d2q17}))
and unknowns $ \, f_j \, $ satisfying a periodicity condition parametrized by a wave      
number $ \, k $: 
\moneq   \label{periodic}
f_j ( x + \xi_j \, \Delta x ,\, t ) \,=\, 
\exp ( i \,\, k \smb \xi_j \,\, \Delta x ) \, \Phi_j  \, , \qquad 0 \leq j \leq q-1 \, . 
\monend 
With a so-called ``zero-point experiment'', 
we consider  one  iteration in time of the d'Hu\-mi\`eres scheme \cite {DdH92} 
with an initial condition satisfying (\ref{periodic}). Such an iteration is  
 defined according to  
\moneq   \label{iteration-temps}
f_j ( x  ,\, t ) \,=\, f_j^* (  x - \xi_j \, \Delta x ,\, t )  
\monend 
with $ \, f_j^* \, $  detailed in Annex~1 at the relation (\ref{f-star}). 
Then the vector $ \, \Phi \, $ is necessarily an eigenvector of the 
amplification matrix, as detailed in \cite{LL00}. The corresponding  eigenvalues 
define the discrete local modes of the linearized scheme. 
They must be  of modulus less than one in order to have a possible stability. 
A typical numerical experiment as those first used in  \cite{LL03}                
is described in Figure~1. For a Prandtl number typically of the order one  
and a wave number greater than a moderate critical wave number,                
the viscous and thermal modes become coupled.                                   
Then the eigenvalues have a non-physical imaginary part, as presented in the
  picture on the right of Figure~1.

\smallskip  \smallskip \noindent 
The physical effects of such bad approximation are  presented {\it e.g.} in figure~2. 
The physical problem is the relaxation of a wave of wave number $ \, k \, $. 
The initial condition is now of the type  
\moneq   \label{k-wave}
f_j ( x  ,\, t ) \,=\, 
\exp ( i \, k \smb x  ) \, \psi_j   \, ,  \qquad  0 \leq j \leq q-1 \,, 
\monend 
with a given vector $ \, \psi \,$ that corresponds to a shear wave and  
 a vertex $ \, x \,$ in a  $ \, N_x \times  N_y \, $ two-dimensional mesh. 
Periodic boundary conditions are enforced. Physically, this wave relaxes towards     
a null equilibrium. For a supercritical wave number, the physics is not well       
approximated by the method: negative values can numerically occur~! 
A consequence of this major default is that very few compressible 
experiments are allowed with the lattice Boltzmann schemes.

\smallskip   \smallskip                        \centerline {  
\includegraphics[width=.52 \textwidth]  {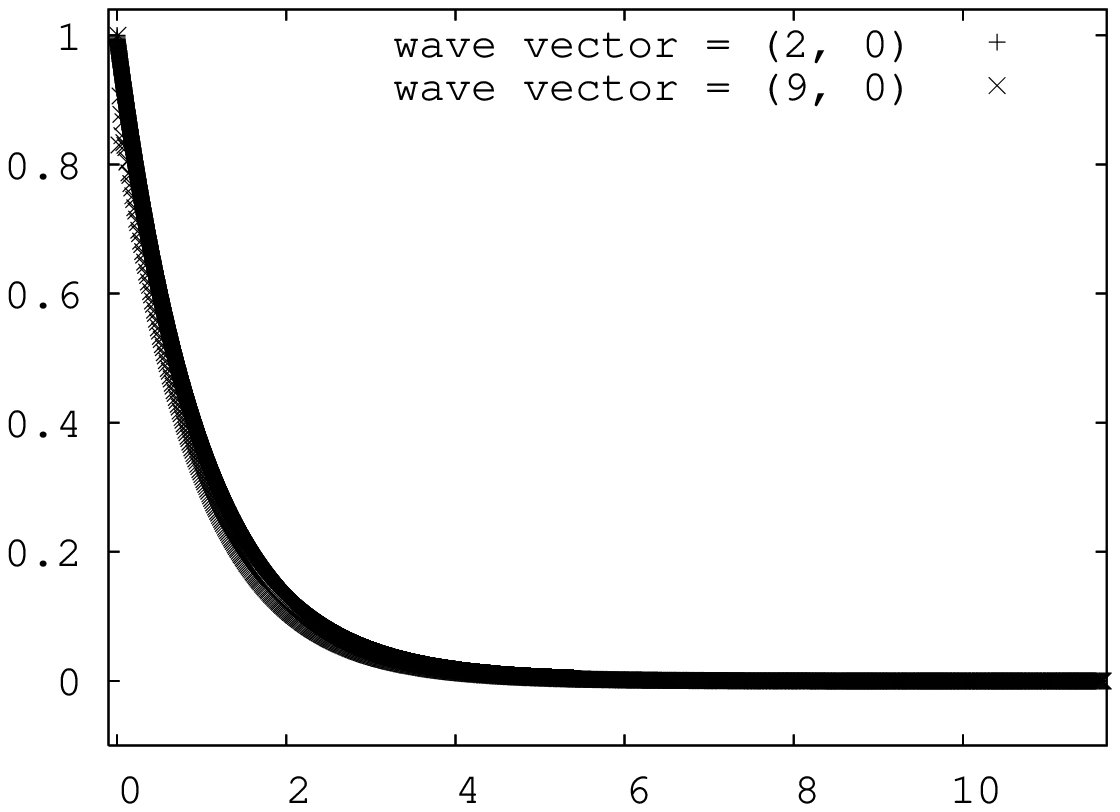} 
\includegraphics[width=.52 \textwidth]  {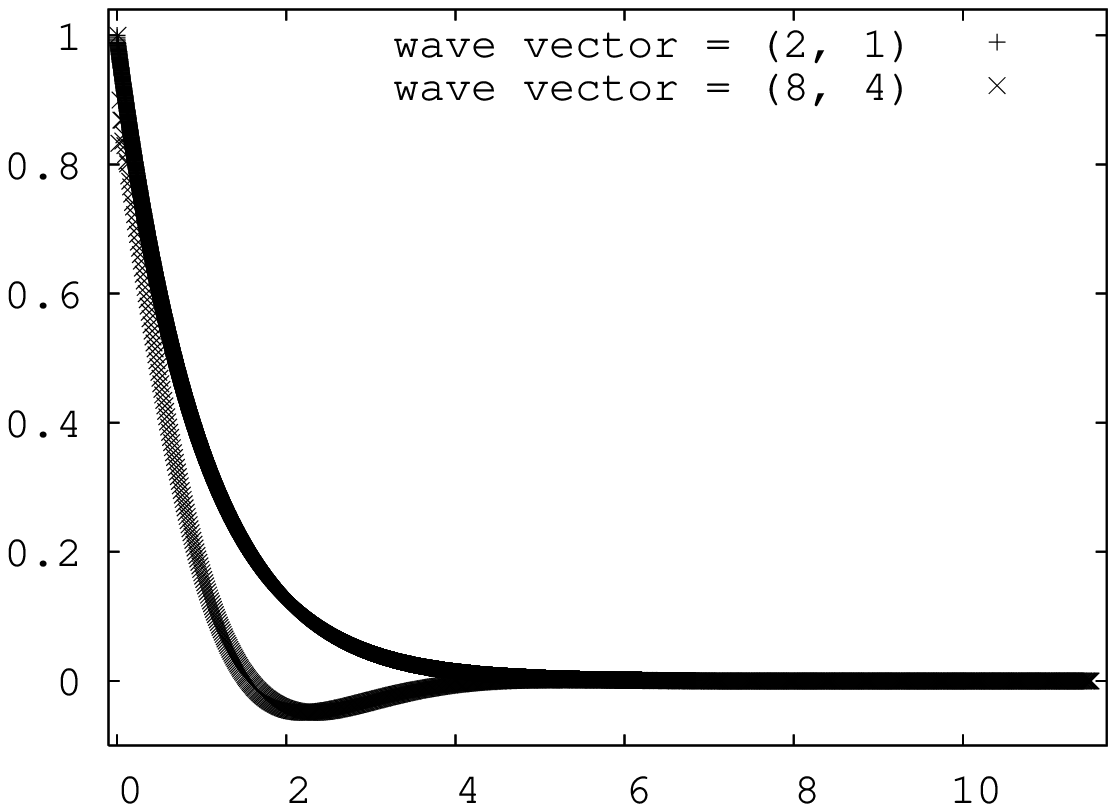} }  

\smallskip \noindent  {\bf Figure 2}. \quad  Relaxation of a thermal wave 
with the D2Q13 lattice Boltzmann scheme.  The amplitude of the wave is plotted as a function   
of the   normalized time. Left: wave vector $\bf k$ parallel to O$x$ axis (2 and 9 wavelengths).  
Right: wave vector $\bf k$ at an angle 26.56$^0$ from O$x$. Two values of the wave number    
 are presented, one smaller than the critical value (see Figure~1) and the other   
larger (respectively 1 and 2 wavelengths along O$x$ and O$y$, and 4 and 8 wavelengths     
along O$x$ and O$y$). The relaxation is physically correct in the first case but an unphysical     
undershoot   appears in the second case. 
Domain $ \, N_x \times  N_y$  with $N_x=N_y=61.$ 
Numerical values of the parameters: 
$ \, s_5=1.88, s_7=1.9303,$\br
 $  s_9=1.60, s_{11}=1.05, s_{12}=1.35, s_{13}=1.65$, 
 $c_1 = -1.3, \alpha_2=-25, \beta_2=-1.5, $\br 
$\alpha_3=4.5, \beta_3= 0.$

\smallskip \smallskip 
  
 \newpage 
\monitem 
In this contribution, we propose some solution to try and solve the previous difficulties. 
We use the Taylor expansion method proposed by one of us \cite{Du08} 
and used in previous contributions for the development of ``quartic'' schemes
  \cite{DL09, DL11}
because the analysis of the full dispersion equation is not practically tractable
when the number $ \, q \,$ of velocities is greater than nine typically.        
With this method, we   analyze the linearized waves of the numerical schemes      
for different problems and the lattice Boltzmann schemes D2Q9, D2Q13 and D2Q17 
presented in Annexes~1, 2 and 3. 
In Section~2,  we show that the inviscid advective acoustics necessarily fixes
some moments of degree 2 and 3.                                      
Then in Section~3,  we consider the second order analysis of the 
lattice Boltzmann scheme  for dissipative  advective
acoustics. We enforce at first order  Galilean invariance 
for shear and thermal waves. 
In Section~4, we analyze the waves of the scheme  
 at fourth order accuracy  for a possible acoustics simulation. 
We enforce isotropy of the waves and this condition fixes an important number      
of parameters  of the method. We propose possible values for the three 
schemes.  In the three following sections, we present  preliminary                  
numerical experiments for the lattice  Boltzmann schemes D2Q9, D2Q13 and D2Q17. 
Some words of conclusion are proposed in Section~8.


 \bigskip \bigskip  \noindent {\bf \large 2) \quad  Inviscid advective acoustics }  

 \noindent    
We are interested by conservation laws of mass, momentum and energy. 
The conserved variables 
\moneq   \label{var-conserv} 
W \,= \, \big( \,\, \rho \,,\,\, j_x \equiv \, \rho \,u   \,,\,\, 
j_y \equiv   \rho \,v  \,,\,\, \varepsilon   \,\, \big) 
^{\displaystyle \rm t}    
\monend 
are related to the particle densities $\, f_j  \,$ through  the relations    
\moneq   \label{moments-conserv}  
\rho \equiv \sum_j \, f_j \,, \,\,\, 
j_x \equiv \sum_j \, v_j^x \, f_j \,, \,\,\, 
j_y \equiv \sum_j \, v_j^y \, f_j \,, \,\,\, 
 \varepsilon \equiv {1\over2} \sum_j \, \mid \! v_j \! \mid ^2 \,  f_j \, + \, 
{\rm orth.} 
\monend 
where ``orth'' are {\it ad hoc} terms for enforcing orthogonality, 
detailed for the various schemes in Annexes 1 to 3. The other second order moments 
are defined by                 
\moneq   \label{xx-xy} 
XX \equiv \sum_j \, \big[ \, (v_j^x)^2 - (v_j^y)^2\,  \big] \, f_j \,, \qquad 
XY \equiv \sum_j \,   v_j^x \, \, v_j^y   \, f_j \, . 
\monend 
The first third order moments $ \, q_x \, $ and  $ \, q_y \, $
are related to heat fluxes: 
\moneq   \label{flux-chaleur} 
q_x \equiv \sum_j \,  {1\over2}  \mid \! v_j \! \mid ^2    v_j^x \, f_j \,, \,\,\,\,\, 
q_y \equiv \sum_j \,  {1\over2}  \mid \! v_j \! \mid ^2    v_j^y \, f_j \, . 
\monend 
In this section these moments at equilibrium are  supposed to be {\bf linearized} functions 
of the conserved moments  $\, W \, $ defined in (\ref{var-conserv}). 
We propose a method for determining   the 16 corresponding coefficients 
when we wish to approximate advective acoustics.

\bigskip \monitem
We start from the Euler equations of gas dynamics 
\moneq   \label{Euler} 
{{\partial W}\over{\partial t}} + {{\partial f(W)}\over{\partial x}} 
 + {{\partial g(W)}\over{\partial y}} \, = \, 0   
\monend 
with 
\moneq   \label{flux-Euler}  \left\{ \begin{array} [c]{l}  
\displaystyle   
f(W) \equiv \big( \, \rho \, u \,, \,\, \rho \, u^2 + p \,, \,\, \rho \, u \, v 
\,, \,\,  u \, \varepsilon  + p \, u \big)^{\displaystyle \rm t}   \\ \displaystyle   
g(W) \equiv \big( \, \rho \, v \,, \,\, \rho \, u \, v  \,, \,\, 
\rho \, v^2 + p \,, \,\,  v \, \varepsilon + p \, v \big)^{\displaystyle \rm t}   \,. 
\end{array} \right. \monend  
We linearize this system around a given state 
\moneq   \label{ref-state}
 W_0     \, = \, \big( \,\, \rho_0 \,,\,\,  \rho_0 \,u_0   \,,\,\, 
  \rho_0 \,v_0  \,,\,\, \rho_0 \, E_0 \,\, \big) ^{\displaystyle \rm t}   . \,   
\monend 
We introduce the internal specific energy $ \, e \,$ according to 
\moneq   \label{energie-interne} 
\varepsilon  \,   \equiv \,  \rho \, e +  {{\rho}\over{2}} \big( u^2 + v^2 \big)    
\monend 
and we suppose that the pressure is  a  function of the  only internal volumic energy    
$ \,  \rho \, e \, $:
\moneq   \label{pression} 
p = {\rm function} \big( \rho \, e \big) \, . 
\monend 
We linearize the pressure $ \, p \,$ given at relation 
(\ref{pression}) around the given state $ \, W_0 \,$ and after 
some lines of elementary calculus, with the notation 
 $ \, \beta \equiv {\rm d}p / {\rm d}(  \rho \, e ) \,$ we have 
\moneq   \label{dp} 
{\rm d}p = \beta_0 \, \Big[  \, {1\over2} \big( u^2 + v^2 \big) \,  {\rm d}\rho 
\,-\, u \, {\rm d}j_x  \,-\, v \, {\rm d}j_y   +  {\rm d}\varepsilon \, \Big] \, . 
\monend 
The three first linearized equations of system (\ref{Euler})(\ref{flux-Euler})     
concerning mass and momentum conservation take the form 
\moneq   \label{Euler-linear}  \left\{ \begin{array} [c]{l} \displaystyle 
 \partial_t \rho  + \partial_x  j_x  +  \partial_y  j_y   = 0    \\ \displaystyle  
 \Big[ \Big(  {1\over2}  \beta_0 ( u_0^2 + v_0^2 ) - u_0^2 \Big) \,  \partial_x 
- u_0 \, v_0 \, \partial_y \Big] \rho  \,+\,  
 \Big[ \partial_t + (2-\beta_0) \, u_0 \, \partial_x +  v_0 \, \partial_y \Big] j_x 
 \\ \displaystyle  \qquad  \qquad  \qquad  \qquad  \qquad 
  \,+\,    \Big[ \! - \beta_0 \, v_0 \,  \partial_x + u_0  \, \partial_y \Big]  j_y 
  \,+\,   \beta_0 \,  \partial_x \varepsilon \,=\, 0   \\ \displaystyle  
 \Big[ \! - u_0 \, v_0 \, \partial_x + 
\Big(  {1\over2}  \beta_0 ( u_0^2 + v_0^2 ) - v_0^2 \Big) \,  \partial_y   \Big] \rho 
 \,+\,     \Big[  v_0  \, \partial_x  - \beta_0 \, u_0 \,  \partial_y  \Big]  j_x 
 \\ \displaystyle  \qquad  \qquad  \qquad  \qquad  \qquad 
  \,+\,  \Big[ \partial_t +  u_0 \, \partial_x +  (2-\beta_0) \, v_0 \, \partial_y \Big] j_y 
  \,+\,   \beta_0 \,  \partial_y \varepsilon \,=\, 0  \, . 
\end{array} \right. \monend  
We identify these equations with those obtained by a first order Taylor expansion   
 (see {\it e.g.} \cite{Du08}) of the lattice Boltzmann scheme. 
Then we obtain for the D2Q9, D2Q13 and  D2Q17 schemes the following expressions 
for second order moments at equilibrium 
\moneq   \label{xx-xy-equil} 
XX^{\rm eq}  = - \big( u_0^2 - v_0^2 \big) \, \rho + 2 \,  u_0 \, j_x - 2 \, v_0 \, j_y  \,, \,\,\, 
XY^{\rm eq}   = -  u_0 \,  v_0 \, \rho  + v_0  \,  j_x     + u_0  \,  j_x     \, . 
\monend 
%

\bigskip \monitem
The  expressions (\ref{xx-xy-equil}) are linear functions of the 
conserved variables  (\ref{var-conserv}) around the 
reference state 
$\, W_0 \,$ given at relation   (\ref{ref-state}). 
If we consider the conserved variables  (\ref{var-conserv}) 
as ``small variations'' of the reference  state  (\ref{ref-state}), 
{\it id est} 
\moneq   \label{infinitesimals} 
 \rho = \delta \rho_0 \,, \,\,   
 j_x  = \delta ( \rho_0 \, u_0) \, , \,\,  
 j_y  = \delta ( \rho_0 \, v_0) \, , \,\,   
 \epsilon  = \delta ( \rho_0 \, E_0) \, , \,\,   
\monend 
and skipping the index ``zero'' for convenience, the  
 expressions (\ref{xx-xy-equil}) can be considered as differential forms:
\moneq   \label{xx-xy-differentielle}   \left\{ \begin{array} [c]{l} \displaystyle 
\delta XX^{\rm eq}  = - \big( u^2 - v^2 \big) \, \delta \rho 
+ 2 \,  u \,  \delta ( \rho \, u) 
 - 2 \, v \,   \delta ( \rho \, v)  \\ \displaystyle 
\delta  XY^{\rm eq}   = -  u \,  v \,  \delta \rho  + v  \,  \delta ( \rho \, u)    
 + u  \,   \delta ( \rho \, v)     \, . 
\end{array} \right. \monend  
A natural question when considering differential forms is to know whether
they are or not the differential of some functions. In other terms, the question 
is to find functions $ \, \xi(\rho,\, u,\, v) \,$ and  $ \, \eta(\rho,\, u,\, v) \,$
such that the expressions given in (\ref{xx-xy-differentielle}) admit also 
the form 
\moneq   \label{xx-xy-diff2}   \left\{ \begin{array} [c]{l} \displaystyle 
\delta XX^{\rm eq}  \equiv \big( u^2 - v^2 \big) \, \delta \rho 
+ 2 \, \rho \,   u \,  \delta  u  -  2 \, \rho \,   v \,  \delta  v 
 \,\, =  \,\,  \delta  \xi(\rho,\, u,\, v)    \\ \displaystyle 
 \displaystyle \delta XY^{\rm eq}   \equiv   u \, v \,  \delta \rho 
+  \rho \,   v \,  \delta  u +  \rho \,   u \,  \delta v
    \,\, =  \,\, \delta  \eta (\rho,\, u,\, v)   \, . 
\end{array} \right. \monend  
If the relations (\ref{xx-xy-diff2}) are true, we have necessarily 
$ \,\,  u^2 - v^2  =  {{\partial \xi}\over{\partial \rho}} \,\, $
and there exists some function $ \, \xi_1(u,\,v) \,$ such that 
$  \,\, \xi (\rho,\, u,\, v) \equiv  \rho \, (  u^2 - v^2 ) +  \xi_1(u,\,v) $. 
Then we have necessarily  
$ \,\,  2 \, \rho \,   u   =  {{\partial \xi}\over{\partial u}} 
=   2 \, \rho \,   u  +  {{\partial \xi_1}\over{\partial u}}  \,\, $
and the function  $ \, \xi_1  \,$ is only function of one
single variable:  $ \, \xi_1 =  \xi_1 (v) $.   
We deduce  from (\ref{xx-xy-diff2}) the new relation 
$ \,\,  - 2 \, \rho \,   v   =  {{\partial \xi}\over{\partial v}} 
=   - 2 \, \rho \,   v  +  {{{\rm d}  \xi_1}\over{{\rm d}  v}}  \,\, $
and  $ \, \xi_1 \,$ is reduced to some constant. 
We can proceed in a similar way for the function  
$ \, \eta (\rho,\, u,\, v) $. First taking the differential 
of the second relation of  (\ref{xx-xy-diff2}) relatively to 
density, we have 
$ \,\,  u \, v   =  {{\partial \eta}\over{\partial \rho}} \,\, $
and there exists some function $ \, \eta_1(u,\,v) \,$ such that 
$  \,\, \eta (\rho,\, u,\, v) \equiv  \rho \,  u \, v   +  \eta_1(u,\,v) $. 
Applying now a derivation relative to $u$: 
$ \,\, \rho \,   v =  {{\partial \eta}\over{\partial u}} 
=   \rho \,  v  +  {{\partial \eta_1}\over{\partial u}}  \,\, $
and $ \,  \eta_1 =  \eta_1 (v) \,$ only. 
After a derivation relative to $v$, we get 
$ \,\, \rho \,   u   =  {{\partial \eta}\over{\partial v}} 
=   \rho \,   u  +  {{{\rm d}  \eta_1}\over{{\rm d}  v}}  \,\, $
and  $ \, \eta_1 \,$ is constant. 
We have proven the relations 
\moneq   \label{xx-xy-diff3}  
\delta XX^{\rm eq}  \, =  \,  \delta  \big( \rho \, (u^2 - v^2 ) \big)  \,, \,\, 
\delta XY^{\rm eq}  \, =  \,  \delta  \big( \rho \, u \, v  \big)  \, . 
\monend  
%
%
The expressions (\ref{xx-xy-diff3}) can be 
integrated up to a constant for nonlinear dynamics (\ref{Euler})(\ref{flux-Euler})
and after a simple change of variables, we obtain nonlinear functions
of the initial conserved variables   (\ref{var-conserv}): 
\moneq   \label{xx-xy-equil-nl} 
XX^{\rm eq}  = {{j_x^2 - j_y^2}\over{\rho}}  \,, \qquad 
XY^{\rm eq}   =  {{j_x \, j_y}\over{\rho}}    \, . 
\monend 
%

\bigskip \monitem
The conservation of energy is more delicate to fit exactly. 
It can be achieved  if we assume that the  equation of state (\ref{pression}) is
precisely  $ \, p =  \rho \, e  \,$ which means that the fluid is a perfect              
gas with a ratio $ \, \gamma \, $ of specific heats equal to 2. In other words, the lattice
Boltzmann schemes are well adapted for shallow water equations. 
For general fluids, we introduce the sound velocity $ \, c_0 \,$ and the Laplace
operator $ \, \Delta \equiv \partial_x^2 + \partial_y^2 . \, $ 
We know that the linearized equations 
\moneq   \label{Euler-linearise}  
A_0 \, \smb \, W \,=\, {\rm O}(\Delta t) 
\monend 
around a given state $ \, W_0 \,$ admit in this  case of two space dimensions 
the following four  eigenvalues
\moneq   \label{val-propres-1}  
\partial_t + u_0 \,\partial_x + v_0 \,\partial_y \,\,\,\, {\rm (double)} , \qquad 
\partial_t + u_0 \,\partial_x + v_0 \,\partial_y \pm c_0 \sqrt{\Delta} \,\,\,\,
  {\rm (acoustics)} 
\monend 
%
with  notations used in \cite{DL11} that are 
 exactly the one used   when implementing the approach with a  symbolic manipulation software. 
It is also possible to introduce a Fourier decomposition on harmonic waves of the type
$\, {\rm exp} \big( i \, (\omega \, t \,-\, {\bf k}  \smb {\bf x} ) \big) . \,$ 
Then we have the usual change of notation: $ \, \partial_t \equiv i \, \omega ,\,$ 
 $ \, \nabla \equiv -i \,  {\bf k} ,\,$  
$ \, \Delta \equiv -  \mid \!\! {\bf k} \!\! \mid^2 ,\,$ 
$ \, \sqrt{\Delta}  \equiv i   \mid \!\! {\bf k} \!\! \mid ,\,$ {\it etc.} 
%

\bigskip \monitem
We impose these eigenvalues to the equivalent equations of the lattice Boltzmann 
schemes D2Q9, D2Q13 and D2Q17. In this way, we obtain 7 independent relationships     
that constrain the equilibrium heat flux  $ \, q \,$ given  at relation (\ref{flux-chaleur})      
 for these three schemes. Independently,  we know   from (\ref{flux-Euler}) 
that when we linearize the conservation of energy, 
the coefficients of $ \, \partial_y  j_x  \,$ 
and  $ \, \partial_x j_y  \,$ are both equal to zero.
In the equivalent equations, we just impose that these two coefficients are equal. 
In this way, we obtain an   
eighth equation.  We solve these equations
and we find for the D2Q9 scheme the following expressions for the linearized heat fluxes
 $ \,  q_x \,$ and  $ \,  q_y \,$ around a given state $ \, W_0 \, $: 
\moneq   \label{flux-linear-D2Q9}  \left\{ \begin{array} [c]{l} \displaystyle 
   q_x^{\rm eq}  \,= 2\, u_0 \, (4 \lambda^2 - 3 c_0^2) \, \rho 
\,+ \, ( 6 c_0^2 + 3 v_0^2 - 3 u_0^2 - 5 \lambda^2 ) \, j_x \,-\, 6 u_0 v_0 \, j_y 
 \,+ \, 2 u_0 \, E   \\ \displaystyle   
   q_y^{\rm eq}  \,= 2\, v_0  \, (4 \lambda^2 - 3 c_0^2) \, \rho
 \,-\, 6 u_0 v_0 \, j_x \,+ \, 
( 6 c_0^2 + 3 u_0^2 - 3 v_0^2 - 5 \lambda^2 ) \, j_y  \,+ \, 2 v_0 \, E  \, . 
\end{array} \right. \monend  
For the D2Q13 lattice Boltzmann scheme, we obtain with the same arguments
\moneq   \label{flux-linear-D2Q13}  \left\{ \begin{array} [c]{l} \displaystyle 
   q_x^{\rm eq}  \,= {{2}\over{13}} u_0 \, (28  \lambda^2 - 13 c_0^2) \, \rho 
\,+ \, ( 2 c_0^2 + v_0^2 -  u_0^2 - 3 \lambda^2 ) \, j_x \,-\, 2  u_0 v_0 \, j_y 
 \,+ \,   {{2}\over{13}} u_0 \, E   \\ \displaystyle   \vspace{-.5cm}  ~  \\ \displaystyle     
   q_y^{\rm eq}  \,=   {{2}\over{13}}  v_0  \,  (28  \lambda^2 - 13 c_0^2) \, \rho
 \,-\, 2 u_0 v_0 \, j_x \,+ \, 
 ( 2 c_0^2 + u_0^2 -  v_0^2 - 3 \lambda^2 )    \, j_y  \,+ \,   {{2}\over{13}}  v_0 \, E    \,   
\end{array} \right. \monend  
and  the D2Q17 scheme leads to 
\moneq   \label{flux-linear-D2Q17}  \left\{ \begin{array} [c]{l} \displaystyle 
   q_x^{\rm eq}  \,= {{6 u_0}\over{17}}  \, (60  \lambda^2 - 17 c_0^2) \, \rho 
- ( 6 c_0^2 + 3 u_0^2 - 3  v_0^2 + 17 \lambda^2 ) \, j_x \,-\, 6  u_0 v_0 \, j_y 
 \,+ \,   {{6 \, u_0 \, E }\over{17}}    \\ \displaystyle   \vspace{-.5cm}  ~  \\ \displaystyle   
   q_y^{\rm eq}  \,=      {{6 v_0 }\over{17}}   \,  (60 \lambda^2 - 17 c_0^2) \, \rho
 \,-\, 6 u_0 v_0 \, j_x + 
 ( 6 c_0^2 + 3 v_0^2 - 3  u_0^2 + 17 \lambda^2 )  \, j_y  \,+\, {{6 \, v_0 \, E }\over{17}} \, .    
\end{array} \right. \monend  
We take into account the relations between the physical total energy 
$ \, \varepsilon \, $ and the orthogonalized  total energy $ \, E \,$ 
presented at relations (\ref{energie-E-epsil-d2q9}), 
 (\ref{energie-E-epsil-d2q13}) and  (\ref{energie-E-epsil-d2q17}). 
%
Using an analysis identical to the one presented in details at the relations 
 (\ref{xx-xy-differentielle}) to (\ref{xx-xy-equil-nl}),  and 
after some lines of elementary calculus, we observe that the relations 
(\ref{flux-linear-D2Q9}), (\ref{flux-linear-D2Q13}) and (\ref{flux-linear-D2Q17})
are linearizations  of the following general relations between the heat flux and the
conserved variables. We have precisely  
\moneq   \label{flux-nonlinear}  \left\{ \begin{array} [c]{l} \displaystyle 
{\rm D2Q9} : \qquad \,\,   {\bf q}^{\rm eq}  \,= \Big(\,   3 \lambda^2 
 \,-\,   3 \mid \! u \! \mid^2   \,+\, 2 \, {{E}\over{\rho}}  \Big) \,  {\bf j}    
 \\ \displaystyle   \vspace{-.5cm}  ~  \\ \displaystyle
{\rm D2Q13} : \qquad  {\bf q}^{\rm eq}  \,= \Big( {{17}\over{13}} \lambda^2 
 \,-\,    \mid \! u \! \mid^2  \,+\,  {{2}\over{13}} {{E}\over{\rho}} \,   \Big) \,  {\bf j}   
 \\ \displaystyle   \vspace{-.5cm}  ~  \\ \displaystyle  
{\rm D2Q17} : \qquad  {\bf q}^{\rm eq}  \,= \Big( {{71}\over{17}} \lambda^2 
 \,-\, 3 \mid \! u \! \mid^2   \,+\, {{6}\over{17}} {{E}\over{\rho}} \,   \Big) \,  {\bf j} \, .  
\end{array} \right. \monend  
We observe at this level of analysis that there is no constraint on the higher order
vectors $ \, {\bf r} \,$ and $\, {\bf \tau} $ whenever they exist
(see the relations (\ref{moments-d2q13}) and (\ref{moments-d2q17})
   of Annexes~2 and~3).

\bigskip \bigskip   \noindent {\bf \large 3) \quad  Dissipative  advective acoustics }  

\noindent 
In the previous section, we have considered the first order eigenvalues given by 
the expressions (\ref{val-propres-1}). We denote by $ \, k_0 \,$ 
the kinetic energy of the reference state:
\moneq   \label{k0} 
k_0 \equiv  {{u_0^2 + v_0^2}\over{2}} \, . \,   
\monend  
Let us set $ \, {\rm u}_0 \smb \nabla \equiv  u_0 \, \partial_x +  v_0 \, \partial_y  \,$ 
and introduce the matrix $ \, \Lambda_0 \, $ as the diagonal 
matrix composed by the eigenvalues:
\moneqstar  
 \Lambda_0 \,=\, {\rm diag} \big( \partial_t +  {\rm u}_0 \smb \nabla \,, \,\, 
\partial_t +  {\rm u}_0 \smb \nabla \,, \,\,   
\partial_t +  {\rm u}_0 \smb \nabla + c_0 \,  \sqrt{\Delta} \,, \,\, 
\partial_t +  {\rm u}_0 \smb \nabla - c_0 \,  \sqrt{\Delta} \, \big) \, . 
\monendstar  
We observe that the corresponding 
matrix of eigenvectors, given according to               
\moneq   \label{R0} 
R_0 \,=\,   \begin{pmatrix} 0 & 1 & \sqrt{\Delta} &  \sqrt{\Delta} \cr
\partial_y &  u_0 & c_0 \,\partial_x + u_0 \,   \sqrt{\Delta}  
& -c_0 \,\partial_x + u_0 \,   \sqrt{\Delta}  \cr
-\partial_x &  v_0 &   c_0 \,\partial_y + v_0 \,   \sqrt{\Delta} &
  -c_0 \,\partial_y + v_0 \,   \sqrt{\Delta}    \cr
u_0 \, \partial_y - v_0 \, \partial_x & k_0  & 
 \displaystyle ( c_0^2 + k_0 ) \,  \sqrt{\Delta} + c_0 \,  {\rm u}_0 \smb \nabla  &
 \displaystyle ( c_0^2 + k_0 ) \,  \sqrt{\Delta} - c_0 \,  {\rm u}_0 \smb \nabla 
\end{pmatrix} \monend  
does {\bf not} depend on the numerical scheme. 
We consider now the equivalent equations of the lattice Boltzmann scheme 
at second order accuracy. With the new variables 
\moneq   \label{var-V} 
V \, = \, R_0 \, \smb \, W 
\monend  
obtained by action of the matrix $ \, R_0 ,\,$  
the equivalent partial differential equations at order~2 take the simple form  
\moneq   \label{edp-ordre-2} 
 \big( \Lambda_0  \,+\, \Delta t \,  P_0  \big) \, \smb \, V 
\,=\, {\rm O}(\Delta t^2)  \, . 
\monend  
The partial differential equations (\ref{edp-ordre-2})                        
 extend naturally the first order expression proposed in (\ref{Euler-linearise}).  
The perturbation matrix $ \, P_0 \,$ 
has the generic form
\moneq   \label{P0} 
P_0 \,=\,   \begin{pmatrix} P_{00}   & \begin{pmatrix} * & * \cr * & * \end{pmatrix} \cr 
\begin{pmatrix} * & * \cr * & * \end{pmatrix}  & 
\begin{pmatrix} * & * \cr * & * \end{pmatrix}   \end{pmatrix} \, . \monend  
The two by two matrix $ \,  P_{00} \, $ 
is {\bf not} diagonal. Then the method of perturbations 
(see {\it e.g.} \cite {Hi91, Ho85})  
that we used in  \cite{DL11} is not straightforward to deal with.      
We have {\it a priori} to diagonalize the perturbation 
 $ \,  P_{00} \, $ which is a difficult task in all generality~!       
In this contribution, following an idea first proposed  by 
Qian \cite{Qi93},   
we want to express that the corresponding 
two first eigenvalues 
\moneq   \label{val-propres-ordre-2} 
\lambda_1 \,=\, 
\partial_t + u_0 \,\partial_x + v_0 \,\partial_y \, + \, \Delta t \, p_1 \,, \quad 
\lambda_2 \,=\, 
\partial_t + u_0 \,\partial_x + v_0 \,\partial_y \, + \, \Delta t \, p_2   \monend  
do not depend on the underlying velocity                               
$ \,  {\rm u}_0 ,\, $ in a way first suggested by Qian and Zhou \cite{QZ98}. 
In this contribution, 
we simply enforce the property that the trace and the determinant 
of the matrix  $ \,  P_{00} \, $ do not depend on $ \,  {\rm u}_0 ,\, $      
at least up to second order. In other terms, we have 
\moneq   \label{condition-vp-ordre-2} 
{{\partial}\over{\partial u_0}} \big( p_j \big) \,=\, 
{{\partial}\over{\partial v_0}} \big( p_j \big) \,=\, 0 \,, \qquad 
j = 1 \,, \, 2 \, .   
\monend  
We did not study the analogous property for acoustic waves, {\it id est} the
condition (\ref{condition-vp-ordre-2}) for $\, j = 3 \, $ and $4$. Such a study will
be considered in future contributions. 
In the end of this section, we explicit the various constraints 
that are obtained for the three lattice Boltzmann schemes due to the conditions 
(\ref{condition-vp-ordre-2}). 

\bigskip   \newpage  \monitem {\bf D2Q9}  

\noindent  
We know first from (\ref {flux-linear-D2Q9}) that there exists some relation between 
the sound velocity $ \, c_0 \, $ and the coefficient $ \, c_1 \,$ defined 
{\it e.g.} thanks to the relation  (\ref{first-equil-D2Q9}). 
The previous relation is enforced and 
the sound velocity is completely imposed:  
\moneq   \label{vson-d2q9} 
c_0 \,=\, \sqrt{2\over3} \, \lambda \, . 
\monend  
%
Moreover, the fifth nonconserved moment 
is the square of energy $ \, E_2 .\,$ It is a scalar field. The conditions 
(\ref{condition-vp-ordre-2}) enforce  this property and we have 
\moneq   \label{E2eq-d2q9} 
E_2^{\rm eq} \,=\, \alpha_2 \, \lambda^4 \, \rho + \beta_2 \, \lambda^2 \, E \, . 
\monend  
Then the perturbations $ \, p_1 \,$ and  $ \, p_2 \,$  define the viscosity and the
diffusivity at constant volume. They are given by             
\moneq   \label{vp-ordre2-d2q9} 
p_1 = - {{\lambda^2}\over{3}}  \, \sigma_5 \, \Delta     \,, \quad 
p_2 = - {{\lambda^2}\over{12}}  \, \big( 4 \,  + 4 \, \beta_2  - \alpha_2   \big) 
\, \sigma_7  \, \Delta \, .  
\monend  
%

\bigskip \monitem {\bf D2Q13 } 

\noindent  
There is {\it a priori} no constraint for the sound velocity. The square of the     
energy at equilibrium is again given by a relation 
of the type (\ref{E2eq-d2q9}). The vectorial moment                              
$ \, {\bf r} \,$ (with labels 8 and 9 in the family (\ref{polynomes-d2q13})) 
is proportional to the momentum $\, {\bf j} \, $:
\moneq   \label{req-d2q13} 
 {\bf r}^{\rm eq} \,=\, {{\lambda^2}\over{12}} \, 
\big( \, 62 \, \lambda^2 - 63 \, c_0^2  \,  \big) \, \,    {\bf j} \, \,
\equiv \,\, c_2 \,  {\bf j} \, . 
\monend  
%
There is no condition for the cube $ \, E_3 \,$ of the energy. The  
13th moment named ``$XX_e$'' is essential for visco-elastic 
simulations when the moment $\, XX \,$ is quasi-conserved. 
It  admits an equilibrium of the type 
\moneq   \label{XXeeq-d2q13}  
XX_e^{\rm eq} \,=\, \xi_x  (u_0 , \, v_0) \,  \, 
\Big(  \, \lambda^4 \, \rho + {{\lambda^2}\over{28}} \, E   \, \Big) \,  . 
\monend  
Remark that we are not completely satisfied by the relation (\ref{XXeeq-d2q13}).
The left and right hand sides are not of the  same tensorial type. 
If we exchange $x$ and $y$, the signs of $XX_e$ is changed but it is not the case
for scalar moments $\rho$ and $\, E .$ Nevertheless, this kind of
lack of tensorial coherence exists at any order if we consider sufficiently 
high order moments, as observed with very different methods by Augier
{\it et al.} \cite {ADG11, ADG12}. 
Then the viscosity and the
diffusivity at constant volume $ \, p_1 \,$ and  $ \, p_2 \,$ take the form
\moneq   \label{vp-ordre2-d2q13} 
p_1 = - {{1}\over{2}}  \, c_0^2 \, \sigma_5 \, \Delta     \,, \quad 
p_2 = - {{1}\over{154}}  \,  {{\lambda^4}\over{c_0^2}} \,
   \big( 28 \, \beta_2  + 140 \,  -  \alpha_2  \big) \, 
\sigma_7  \, \Delta \, .  
\monend  
%

\bigskip \monitem {\bf D2Q17 }

\noindent  
As for the D2Q13 scheme, there is no  constraint for  the sound velocity.    
The square of the  energy at equilibrium is still obtained  by the condition  
 (\ref{E2eq-d2q9}). 
There is no condition for the ``powers'' three  $ \, E_3 \,$ and 
four   $ \, E_4 \,$  of the energy.
Note that  ``$XX_e$'' and  ``$XY_e$'' (labels 12 and  13  in (\ref{polynomes-d2q17})) 
satisfy conditions close to (\ref{XXeeq-d2q13}):
\moneq   \label{XXeeq-d2q17}  
XX_e^{\rm eq} \,=\,   \xi_x  (u_0 , \, v_0)  \,  \, 
\Big(  \, \lambda^4 \, \rho + {{\lambda^2}\over{60}} \, E   \, \Big) \, , \quad 
XY_e^{\rm eq} \,=\,    \xi_y  (u_0 , \, v_0)  \,  \, 
\Big(  \, \lambda^4 \, \rho + {{\lambda^2}\over{60}} \, E   \, \Big) \, . 
\monend  
There is no condition on the equilibrium of vector $ \, {\bf r} .\,$ 
But if we introduce the notations 
\moneq   \label{req-d2q17} 
 \left\{   \begin{array} [l]{l}  \displaystyle 
 r_x^{\rm eq} \,=\, \lambda^3 \, \big( c_x^\rho  \, \lambda^2   \, \rho +   c_x^x \, \lambda \, j_x  
 +   c_x^y  \, \lambda \, j_y  +     c_x^\varepsilon   \, E \big)    \\ \displaystyle 
r_y^{\rm eq} \,=\, \lambda^3 \, \big(  c_y^\rho  \, \lambda^2   \, \rho +   c_y^x     \, \lambda \, j_x  
+     c_y^y \, \lambda \, j_y  +     c_y^\varepsilon  \, E \big)  \,,    
\end{array} \right.     \monend
where the $\, c's \,$ coefficients of relations (\ref{req-d2q17}) are 
{\it a priori} functions of the advection field $ \, {\bf u_0} $, 
we have the following expressions for the vector                        
$ \, {\bf \tau} = (  \tau_x  \,,\,   \tau_y  )   \equiv 
( X \, E_3 \, + {\rm orth}. \,,\,  Y \, E_3  \, + {\rm orth}. ) \,$ 
with labels 10 and 11 at relations  (\ref{polynomes-d2q17}): 
\moneq   \label{taueq-d2q17} 
 \left\{   \begin{array} [l]{l}  \displaystyle 
\tau_x^{\rm eq} \,=\, -  {{31}\over{2}} \, \lambda^5 \, 
\Big[   c_x^\rho \, \lambda^2   \, \rho +    {{\lambda}\over{124}} 
\Big( 124 \,   c_x^x + 249 \, {{c_0^2}\over{\lambda^2}} - 442  \Big)  \, j_x  
+   c_x^y \, \lambda \, j_y  +   c_x^\varepsilon  \, E \, \Big]   
 \\ \displaystyle   \vspace{-.5cm}  ~  \\ \displaystyle
\tau_y^{\rm eq} \,=\,  -  {{31}\over{2}} \, \lambda^5 \, 
\Big[    c_y^\rho  \, \lambda^2    \, \rho +      c_y^x  \, \lambda \, j_x  + 
  {{\lambda}\over{124}} \Big( 124 \,   c_y^y + 249 \, {{c_0^2}\over{\lambda^2}} 
- 442  \Big)   \, j_y             +    c_y^\varepsilon \, E \, \Big]    \,.    
\end{array} \right.     \monend
Finally  the perturbations $ \, p_1 \,$ and  $ \, p_2 \,$ are given by     
\moneq   \label{vp-ordre2-d2q17} 
p_1 = - {{1}\over{2}}  \, c_0^2 \, \sigma_5 \, \Delta     \,, \quad 
p_2 = - {{1}\over{218}}  \,  {{\lambda^4}\over{c_0^2}} \,
   \big( 60 \, \beta_2  + 620 \,    -  \alpha_2   \big) \, \sigma_7  \, \Delta \, .  
\monend  
A variant of the relations (\ref{vp-ordre2-d2q13})~!

\bigskip \bigskip   \noindent {\bf \large 4) \quad  Fourth order isotropic  acoustics }  

\noindent 
We suppose in this section that the reference advective 
state $ \, W_0 \,$ has a zero velocity~: $ \, u_0 = v_0 = 0 $. 
We evaluate  the eigenvalues $ \, \lambda_j \,$ (for $j = 1$ to 4) 
at fourth order accuracy by using the general method
presented in details in \cite{DL11}. 
Then the eigenvalues admit a general expansion of the type 
\moneq   \label{vp-ordre4} 
 \left\{   \begin{array} [l]{l}  \displaystyle 
\lambda_1 \,=\, \partial_t + \Delta t \, p_1  + \Delta t^2  \, \widetilde{p_1} 
 + \Delta t^3  \, \overline{p_1} + {\rm O}( \Delta t^4 ) \\ \displaystyle   
\lambda_2 \,=\, \partial_t + \Delta t \, p_2  + \Delta t^2  \, \widetilde{p_2} 
 + \Delta t^3  \, \overline{p_2} + {\rm O}( \Delta t^4 ) \\ \displaystyle   
\lambda_3 \,=\, \partial_t + c_0 \, \sqrt{\Delta} 
+ \Delta t \, p_3  + \Delta t^2  \, \widetilde{p_3} 
 + \Delta t^3  \, \overline{p_3} + {\rm O}( \Delta t^4 ) \\ \displaystyle   
\lambda_4 \,=\, \partial_t - c_0 \, \sqrt{\Delta} 
+ \Delta t \, p_3  - \Delta t^2  \, \widetilde{p_3} 
 + \Delta t^3  \, \overline{p_3} + {\rm O}( \Delta t^4 ) 
\end{array} \right.     \monend
and we refer to 
(\ref{val-propres-1}) and (\ref{val-propres-ordre-2}) 
for advective acoustics at first and second order accuracy. 
In the following, we enforce isotropy by saying that the eigenvalues
$ \, \lambda_j \,$ proposed in (\ref{vp-ordre4}) are isotropic.
In other words, the operators $ \, p_j $ , $\,  \widetilde{p_j} \,$ 
and $ \,  \overline{p_j} \,$ that appear in  (\ref{vp-ordre4})
are only functions of the Laplacian. This induces a family of equations 
for the parameters.

\bigskip \monitem {\bf D2Q9 } at third order accuracy. 

\noindent  For the D2Q9 lattice Boltzmann scheme, we have a total of 
5~equations (respectively one equation) to achieve isotropy at the fourth 
(respectively   third) order.  
We have no solution at the fourth order. Third order isotropy can be enforced, 
{\it i.e.} the dispersion of ultrasonic waves is isotropic  in this case, 
by adding to the relations (\ref{vson-d2q9}) and  (\ref{E2eq-d2q9})
the constraint 
\moneq   \label{sisi-s12} 
\sigma_7 = {1\over{12\, \sigma_5}} \, . \, 
\monend  
%

\bigskip  \newpage \monitem {\bf D2Q13 } 

\noindent For this scheme, a total of  6 equations is necessary to obtain fourth order
isotropy. They can be solved by adding to the previous conditions 
(\ref{req-d2q13}) and (\ref{XXeeq-d2q13}) 
the constraint (\ref{sisi-s12}) and the following specific relations 
\moneq   \label{vson-d2q13}   
c_0 \,=\, {{2}\over{\sqrt{5}}} \, \lambda     \, , \qquad 
\monend  
and
\moneq   \label{E3-d2q13}  
E_3^{\rm eq}  \, = \,  \alpha_3 \, \lambda^6 \, \rho + \beta_3 \, \lambda^4 \, E   \, .
\monend  
The coefficients  $ \,  \alpha_3 \,$ and $ \,  \beta_3 \,$ of the relation
(\ref{E3-d2q13}) are associated to  the coefficients 
 $ \,  \alpha_2 \,$ and $ \,  \beta_2 \,$ introduced at relation (\ref{E2eq-d2q9})  according to
\moneq   \label{alpha_3-beta_3-d2q13} 
 \left\{   \begin{array} [l]{l}  \displaystyle 
\alpha_3 = {{1}\over{1716}} \,\, {{N_\alpha } \over{384 \, \sigma_5^2 + 7}}   \,\,,  \qquad 
\beta_3 =   {{1}\over{216216}} \,\,  {{ N_\beta } \over{ 384 \, \sigma_5^2 + 7  }}   \,\,, 
 \\ \displaystyle   \vspace{-.5cm}  ~  \\ \displaystyle   
N_\alpha = 41922 - 2505 \, \alpha_2 + 54800\, \beta_2 
+ \big(  14098944 \,  + 97440 \, \alpha_2 + 1315200 \, \beta_2  \big) \, \sigma_5^2  
 \\ \displaystyle   \vspace{-.5cm}  ~  \\ \displaystyle   
N_\beta =  - 2756851  + 34250 \, \alpha_2 -889970 \, \beta_2  
\\ \displaystyle \qquad \qquad \qquad \qquad \qquad \qquad 
-  \big( 204329472  - 822000 \, \alpha_2  +  41211840 \, \beta_2  \big) \,  \sigma_5^2 
\,  .  
\end{array} \right.     \monend
The coefficient $ \, \xi_x \,$ in the relation (\ref{XXeeq-d2q13})   is null and we have also  
the following relations between the dissipation coefficients defined  in (\ref{sigma-s})
from the $ \, s's \,$: 
\moneq   \label{sigma9-13}  
\sigma_9 = \sigma_7 \,, \qquad 
\sigma_{13} = \sigma_5    \, 
\monend  
%

\bigskip \monitem {\bf D2Q17 }

\noindent  In this case, fourth order isotropy induces a total of 9~equations. 
They can be solved analytically (with the help of a formal software for the algebra)
first by considering the relations  (\ref{E2eq-d2q9}) and  (\ref{XXeeq-d2q17}).
Secondly, the sound velocity $ \, c_0 \,$ has not to be imposed. 
We have to enforce (\ref{sigma9-13}) and we add  the condition 
\moneq   \label{sigma15-d2q17}   
\sigma_{15} = \sigma_7    \,.   
\monend  
Relation (\ref{E3-d2q13}) is supplemented by an analogous one for the fourth power of the 
energy:
\moneq   \label{E4-d2q17}  
E_4^{\rm eq}  \, = \,  \alpha_4 \, \lambda^8 \, \rho + \beta_4 \, \lambda^6 \,E   \,. \,  
\monend  
The coefficients $ \alpha$'s and  $ \beta$'s satisfy now 
\moneq   \label{alphabet-d2q17} 
 \left\{   \begin{array} [l]{l}  \displaystyle 
\alpha_3 = -  {{5}\over{436}} \,  (2696442 \,   + 7519 \,   \alpha_2 )   \,, \quad  
\beta_3 = -  {{1}\over{2616}} \,  (2949247 \,  + 225570 \,  \beta_2 )  \,, 
 \\ \displaystyle   \vspace{-.5cm}  ~  \\ \displaystyle
\alpha_4 = -  {{1}\over{177888}} \,  (69687842   + 139145 \,   \alpha_2 )   \,, \,\,   
\beta_4 = -  {{5}\over{355776}} \,  (940101    + 55658 \,   \beta_2) \, .  
 \end{array} \right.  
\monend  
Moreover, the vectors $ \,  {\bf q} ,\,$  $ \, {\bf r} \,$ and 
$ \, \tau \,$ considered previously satisfy at equilibrium               
the relations 
\moneq   \label{vecteurs-d2q17}  
 \left\{   \begin{array} [c]{c}  \displaystyle 
 {\bf q}^{\rm eq}  \,=\, c_1 \, {\bf j} \,, \qquad 
 {\bf r}^{\rm eq}  \,=\, c_2 \, {\bf j} \,, \qquad 
 \tau^{\rm eq}   \,=\, c_3 \, {\bf j}  \\ \displaystyle    
c_1 = 6 \, c_0^2 -17 \, \lambda^2 \,, \quad 
c_2 = {{\lambda^2}\over{6}} \, ( 31 \, \lambda^2 - 21 \, c_0^2 )  \,, \quad 
c_3 = {{\lambda^4}\over{24}} \, ( 555 \,  c_0^2  - 596 \, \lambda^2 )   \,.
  \end{array} \right.  \monend   
We have also the simple
equilibria 
\moneq   \label{xxe-xye-d2q17}  
 XX_e^{\rm eq} \,=\, 0 \,, \qquad  XY_e^{\rm eq}  \,=\, 0 \, . 
\monend   
%

\bigskip \bigskip   \noindent {\bf \large 5) \quad  Numerical experiments with the D2Q9 scheme} 

\noindent 
We first consider a ``zero point'' analysis as described in the introduction.
We observe in Figure~3 (left)                                          
that the unphysical coupling of waves is present with an arbitrary 
value of the parameter $ \, \sigma_7 \,$ which is proportional to the diffusivity 
$ \, \kappa \, $ at constant volume 
as indicated at the relation (\ref {vp-ordre2-d2q9}).                               
When fourth order isotropy is enforced according to the relation  (\ref{sisi-s12}), 
this coupling disappears, as observed in Figure~3 (right).               

\smallskip   \smallskip   \centerline {                      
\includegraphics[width=.50 \textwidth]  {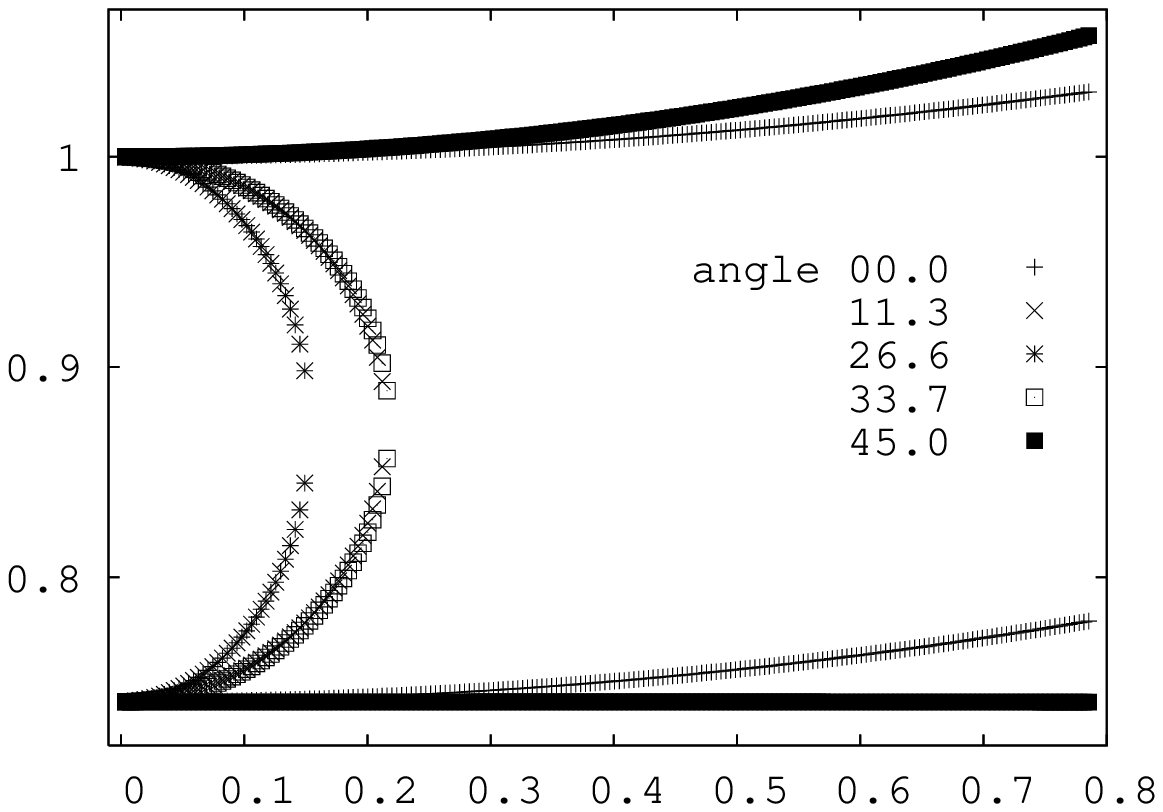} 
\includegraphics[width=.50 \textwidth]  {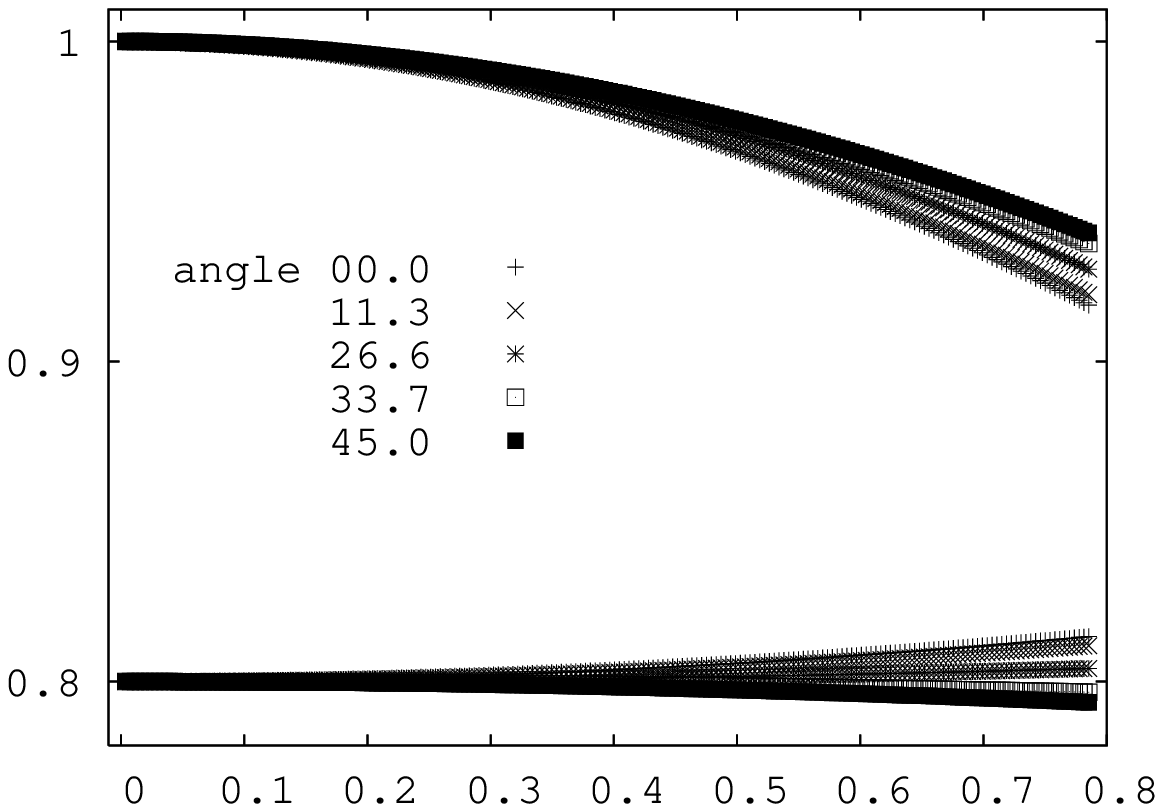} }

\smallskip \noindent  {\bf Figure 3}. \quad D2Q9 ``zero point''. 
Value of the eigenmode divided by $ \, k^2 \,$ 
and normalized by the diffusivity   $ \kappa $ {\it vs} the wave number. \quad      
Left figure~: shear and thermal waves 
with $ \, \sigma_7 \,$   chosen arbitrarily. We see clearly a strong coupling between the   
viscous and diffusive waves for an angle $ \, \theta = 26.565 \,$ degrees. 
Right figure~: the relation (\ref{sisi-s12})  is satisfied.   
The coupling has disappeared but there is still an angular dependency
that characterizes this third order isotropy.       
\smallskip \smallskip 

\smallskip   \smallskip  \centerline  {  
$ \!\!\!\!\!\! $
\includegraphics[width=.33 \textwidth] {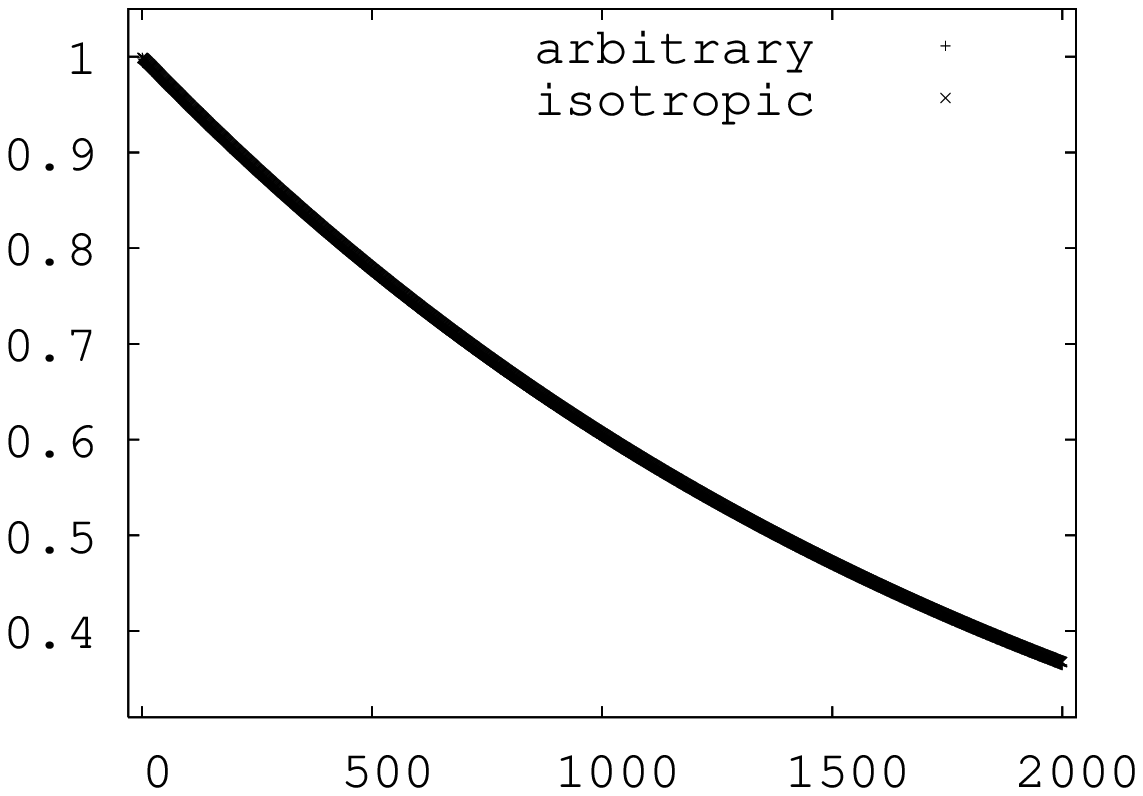}  
\includegraphics[width=.33 \textwidth] {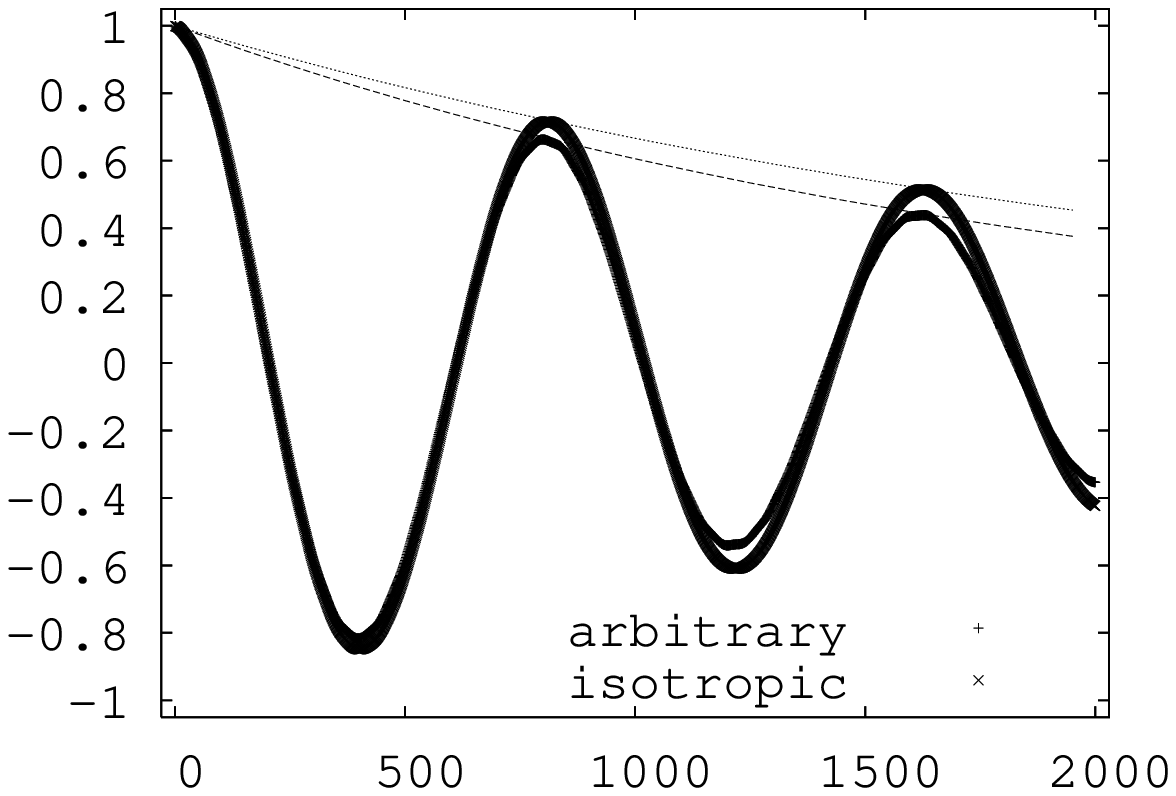}  
\includegraphics[width=.33 \textwidth] {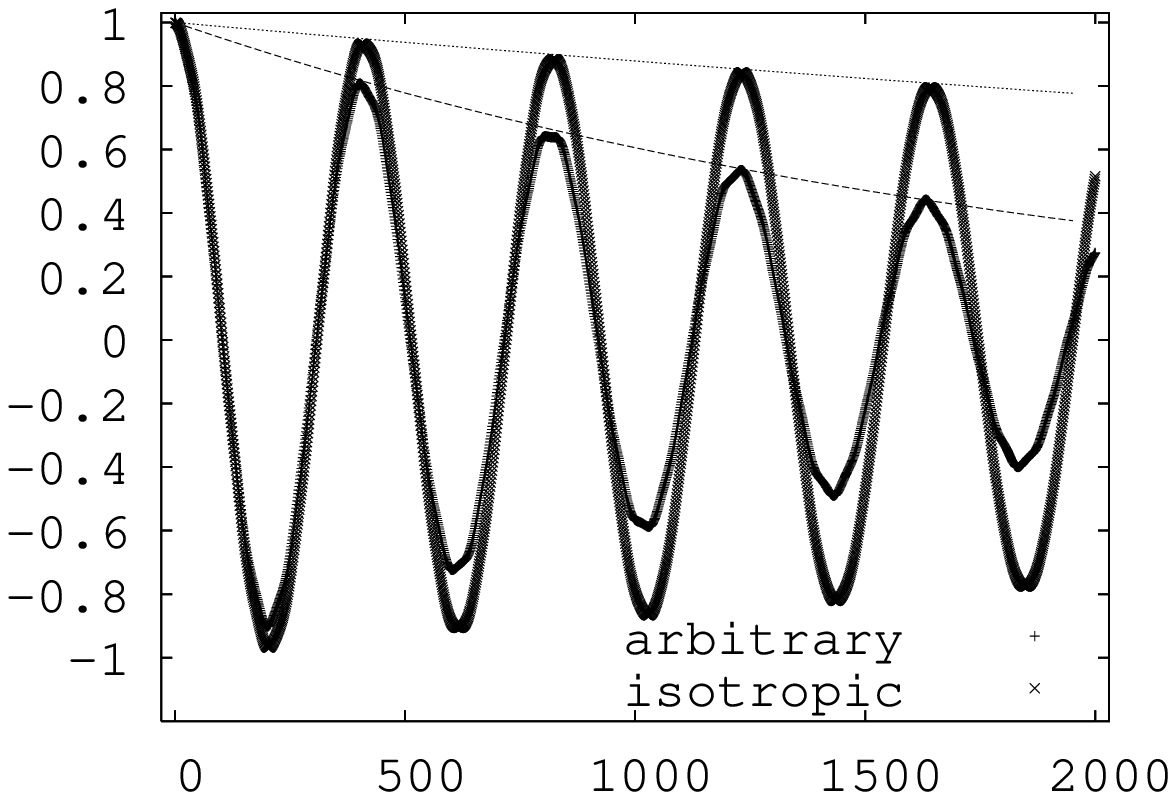} }

\smallskip \noindent  {\bf Figure 4}. \quad Relaxation of a nonlinear diffusion wave with 
D2Q9 ``energy conserving'' lattice Boltzmann scheme.
The parameters of the scheme are the following: 
$ \, s_5 = 1.8181 ,\,$   $ \, s_9 = 1.1765,\,$ 
$ \, \alpha_2 = -1 ,\,$ $ \, \beta_2 = 0.1 \,$ (see the relation (\ref{E2eq-d2q9})), 
$ \, s_7 = 0.4615 \,$ when the condition  (\ref{sisi-s12}) is not satisfied,  
$ \, \alpha_2 = -0.15 ,\,$ $ \, \beta_2 = -1 \,$ (see the relation (\ref{E2eq-d2q9})), 
$ \, s_7 = 1.8305  \,$ when the condition  (\ref{sisi-s12}) is satisfied.
Left:    $v_{advec} = 0 , \,  $ 
Middle:    $v_{advec} = 0.05 $, 
Right:    $v_{advec} = 0.10 $.       
The light exponential lines correspond to the velocity corrected damping (following        
complicated expressions not given here).              
 \smallskip \smallskip 

\monitem 
In order to confirm this good performance of the D2Q9 lattice Boltzmann scheme with
conservation of energy, we have simulated the relaxation of a thermic wave on a     
81 $ \times$ 81 lattice.                                       
We have incorporated the   nonlinear terms given by      
 relations (\ref{xx-xy-equil-nl}) for the moments $ \, XX \, $ 
and  $ \, XY \, $ at equilibrium. 
For the  ``heat flux'' $ \, {\bf q} \,$ at equilibrium, 
we have considered the expressions (\ref{flux-nonlinear}), but  
 the quadratic term relative to velocity has been neglected. 
The results are presented in Figure~4. 
For a small wave vector $k$ and an advection velocity $V$
parallel to the wave vector, 
the  waves are correctly advected,                               
whatever the direction of the wave vector. 
In other terms, we have  isotropy  of the Galilean factor.

\bigskip \bigskip  \bigskip  
 \noindent {\bf \large 6) \quad  Numerical experiments with the D2Q13 scheme}
 

\noindent 
With the methodology presented in sections 2 to 4, it is possible to
remove the spurious coupling of shear and thermal modes depicted in the introduction.   
Precisely, if the parameters of the scheme 
satisfy the relations    (\ref{req-d2q13}), (\ref{sisi-s12}),   (\ref{vson-d2q13}), 
    (\ref{E3-d2q13}),  (\ref{sigma9-13}),  (\ref{c1-d2q13}), 
there exists a situation where the scheme is linearly stable for     
fluid and thermal applications and also for pure acoustics.
Moreover we obtain a  correct Prandtl number and 
appropriate attenuations : 
\moneqstar  
Pr = 0.728 \,, \quad 
\nu  = 0.006  \, \lambda \, \Delta x  \,, \quad 
\kappa = 0.008236  \, \lambda \, \Delta x  \,, \quad 
\gamma = 0.003487   \, \lambda \, \Delta x  \quad          
\monendstar
The results are proposed in Figure~5.                      

\smallskip    \bigskip  
\centerline { \includegraphics[width=.55  \textwidth] {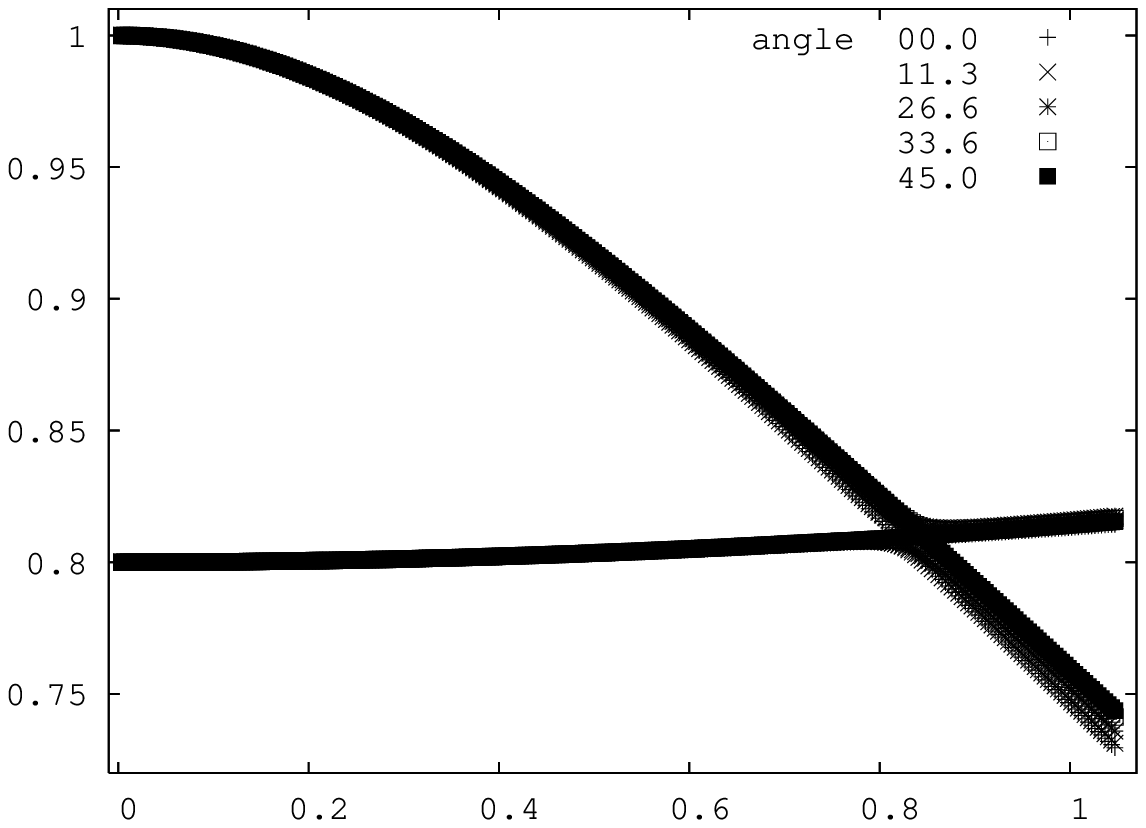}}

\smallskip \noindent  {\bf Figure 5}. \quad  ``Zero point'' experiment 
with the D2Q13 lattice Boltzmann scheme.  
Viscous and diffusive modes for a moderate wave number $ k $  and several angles.     
There is clearly isotropy and the two waves are decoupled.  
The diffusive wave at $k=0$ is on the order of 1.025.  Note that the small oscillations
at this point reflect the numerical difficulties due to the approximation
of the eigenvalue $1$ at fourth order accuracy.  
The viscous wave  at $k=0$ is on the order of 0.72. 
At  $k \simeq 0.78 $ the two modes cross perfectly without merging. 
Note that this perfectly isotropic test case is also very dispersive. 
 
\newpage 

~

\monitem   The relaxation of a diffusive wave is presented Figure~6. 
 
\vskip -.5cm 

  \smallskip     \smallskip   

 \centerline { 
\includegraphics[width=.33 \textwidth] {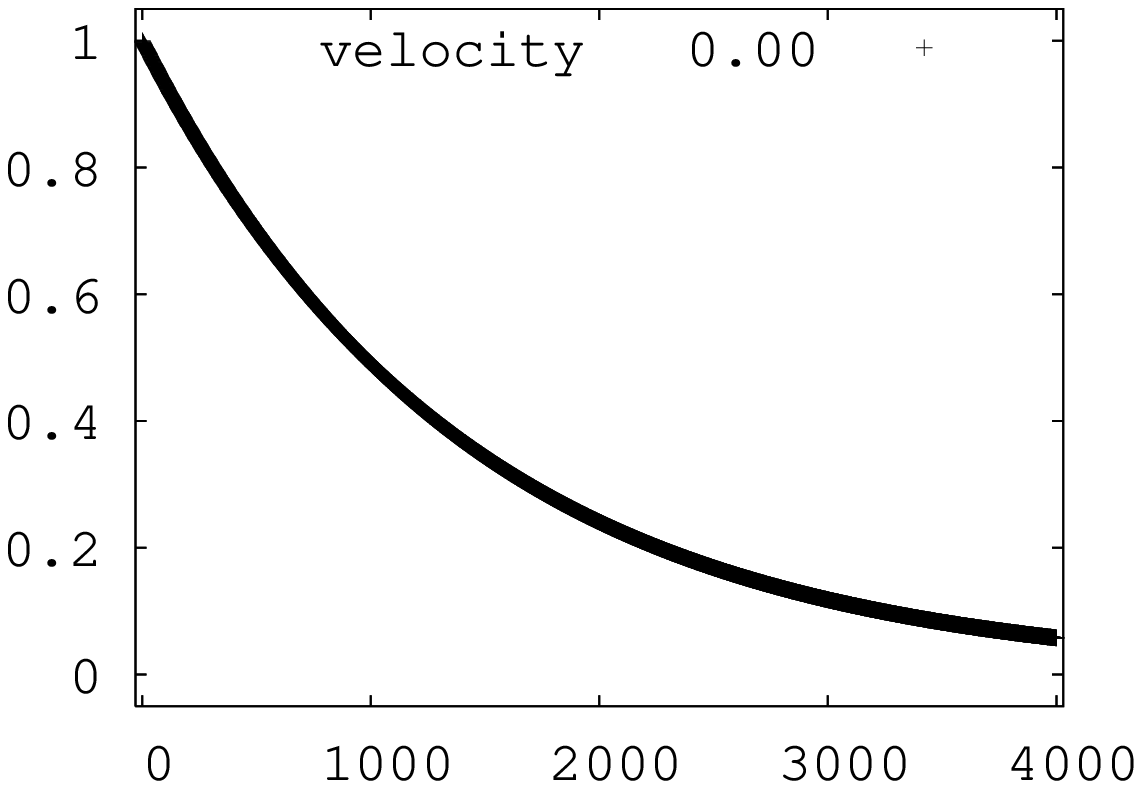} 
\includegraphics[width=.33 \textwidth] {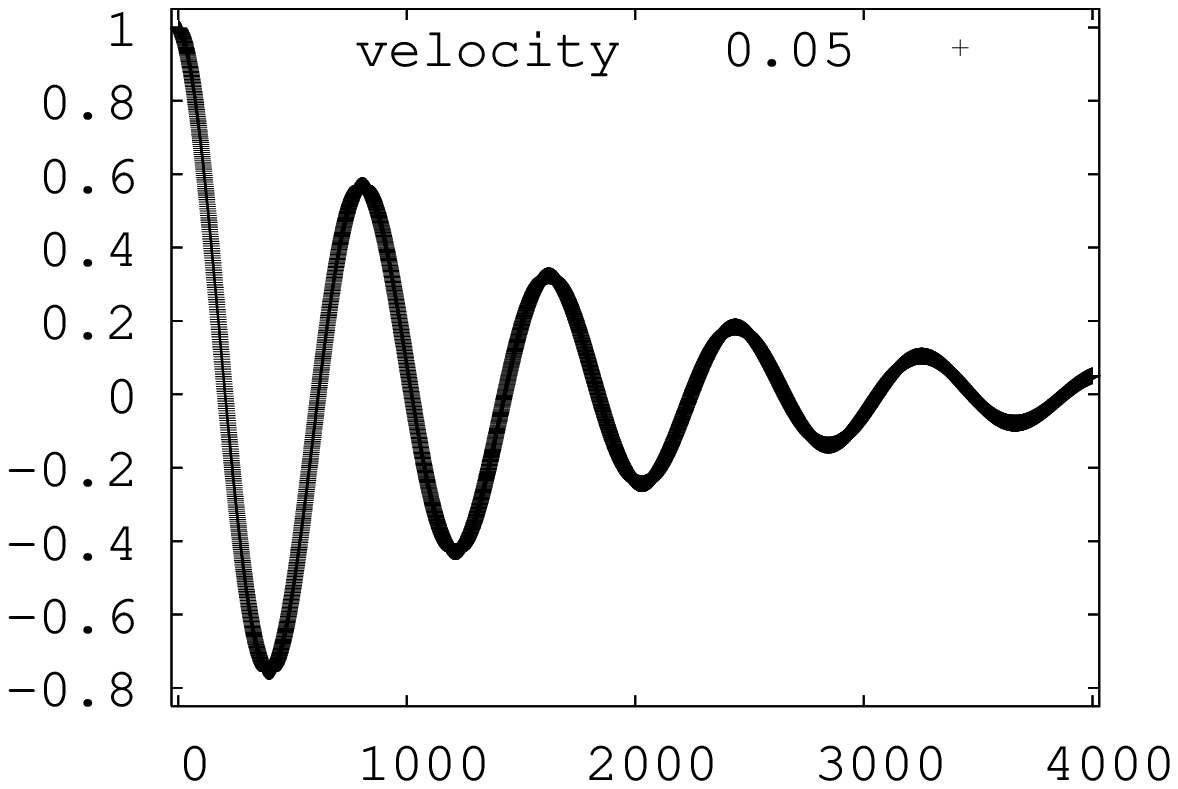}  
\includegraphics[width=.33 \textwidth] {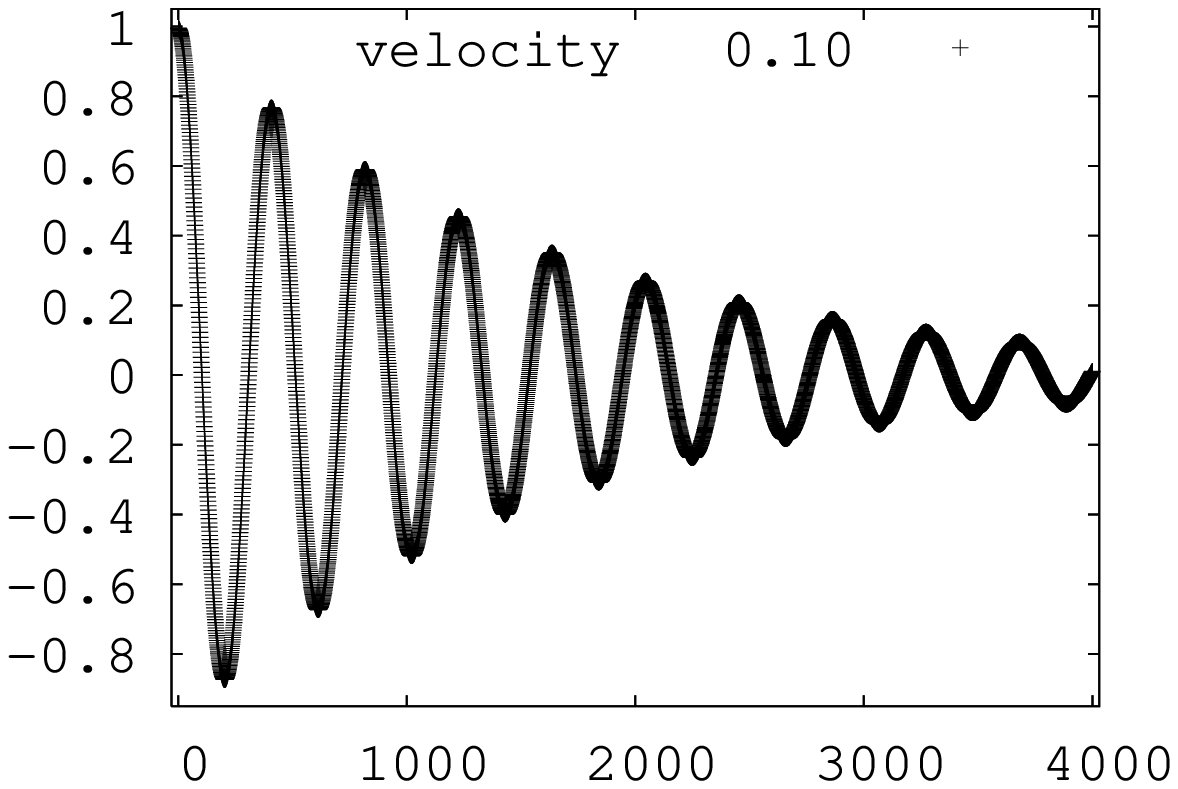} $\,$ }

 \smallskip  \smallskip  \noindent  {\bf Figure 6}. 
\quad Relaxation of a nonlinear diffusive wave as function of time with the 
D2Q13 ``energy conserving'' lattice Boltzmann scheme in a 91 $\times$ 91 domain (2 wave lengths along
O$x$ and 1 wave length along O$y$).
Parameters are set to have a Prandtl number of .80. Mean velocity parallel to the wavevector
of amplitude 0.0, 0.05 and 0.10.

 \smallskip \smallskip \bigskip 

\monitem 
As an illustration of the potential of this ``conserving energy lattice Boltzmann
scheme'', we  present in Figure~7  the propagation of a sound wave in a disc.         

\smallskip   \smallskip   \bigskip                         
  \centerline { \includegraphics [width=.45 \textwidth,] {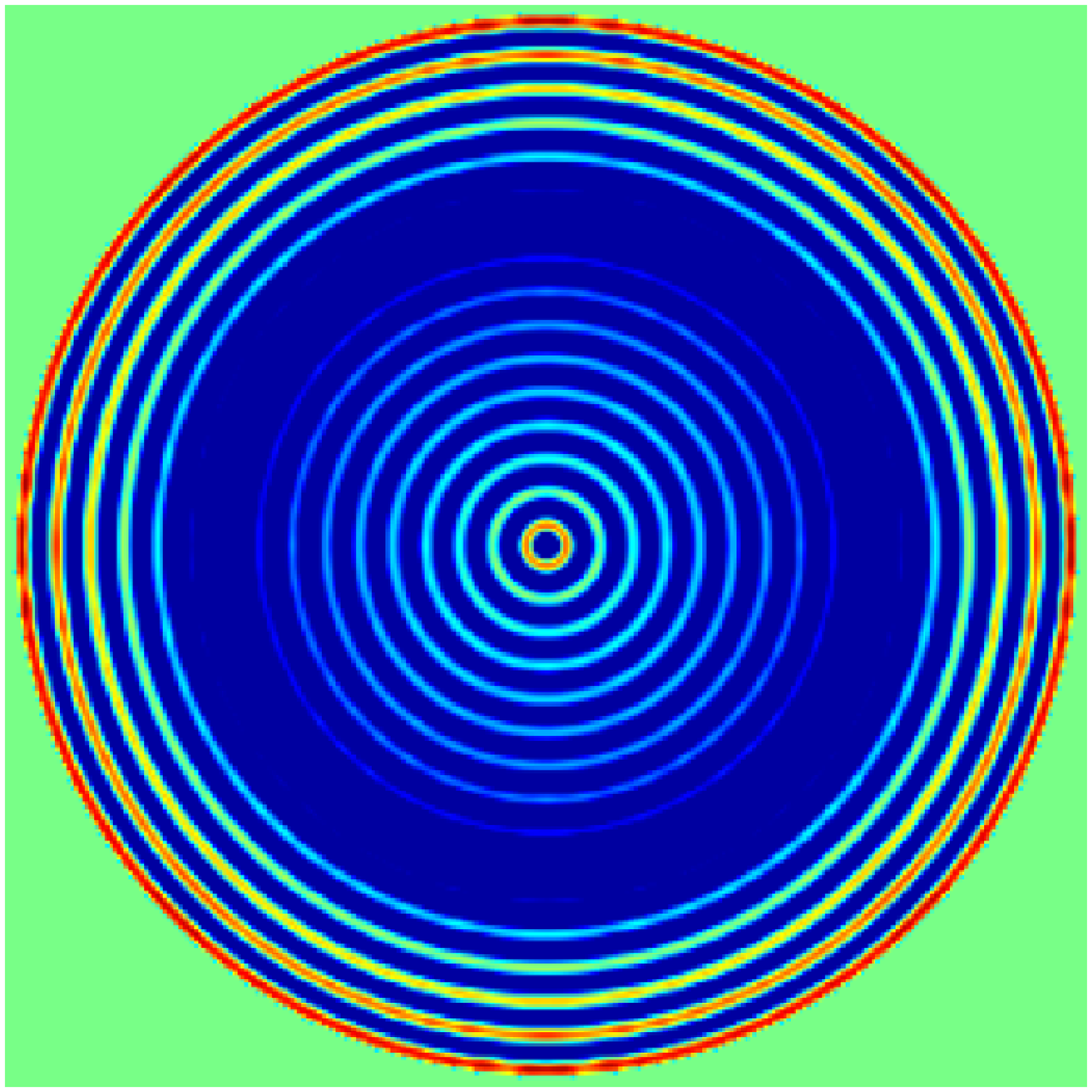}  }

\smallskip \noindent  {\bf Figure 7}. \quad   
Sound wave propagation in a circle with an ``anti bounce - back'' numerical boundary
condition with the D2Q13 lattice Boltzmann scheme conserving the energy.   

\bigskip \bigskip   \noindent {\bf \large 7) \quad  Numerical experiments with the D2Q17 scheme}

\noindent
With the methodology presented in Sections~2 to 4, 
the D2Q17 scheme depicted in Figure~10 and in Annex~3 admits parameters      
satisfying the numerical constraints made explicit  in relations           
(\ref{xx-xy-equil}), 
(\ref{xx-xy-equil-nl}), 
(\ref{flux-nonlinear}),
(\ref{XXeeq-d2q17}), 
(\ref{req-d2q17}), 
(\ref{taueq-d2q17}),
(\ref{sigma15-d2q17}), 
(\ref{E4-d2q17}),
(\ref{alphabet-d2q17}), 
(\ref{vecteurs-d2q17}), 
(\ref{xxe-xye-d2q17}) and 
(\ref{sigma-d2q17}).

\bigskip

\bigskip \monitem 
Some results are shown for the ``zero-point'' analysis.                     
In a first case, we have taken the parameters in a simple way.      
The results are presented in Figure~8. 
The shear and eigenmodes are decoupled and show very little angular dependence      
%
The decoupling of viscous and thermal modes is correct.                        
 %
\moneqstar  
Pr = 0.74182  \,, \,\, 
c_0 =  \sqrt{7\over6}  \,, \,\, 
{{\nu}\over{\lambda \, \Delta x}}   =   0.029167       \,, \,\,  
{{\kappa}\over{\lambda \, \Delta x}} =   0.039318   \,, \,\,  
{{\gamma}\over{\lambda \, \Delta x}} =   0.055959 \, .      
\monendstar

\smallskip   \smallskip              
 \centerline {  
\includegraphics[width=.33 \textwidth] {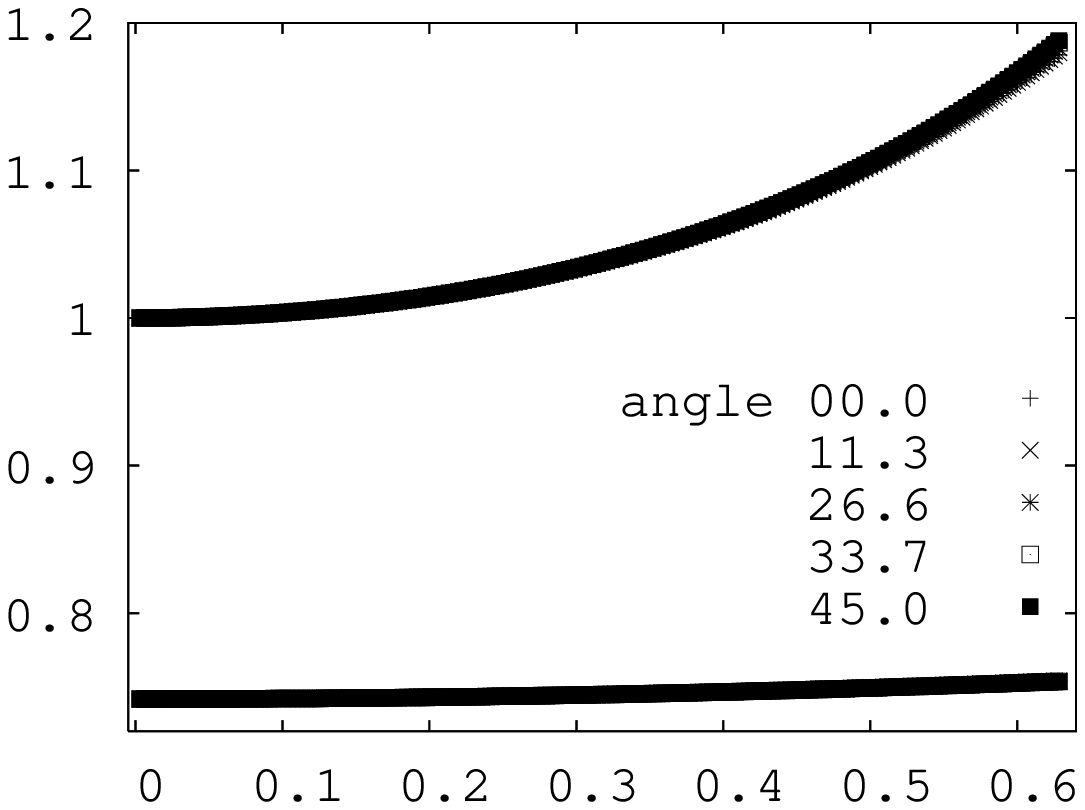}  
\includegraphics[width=.33 \textwidth] {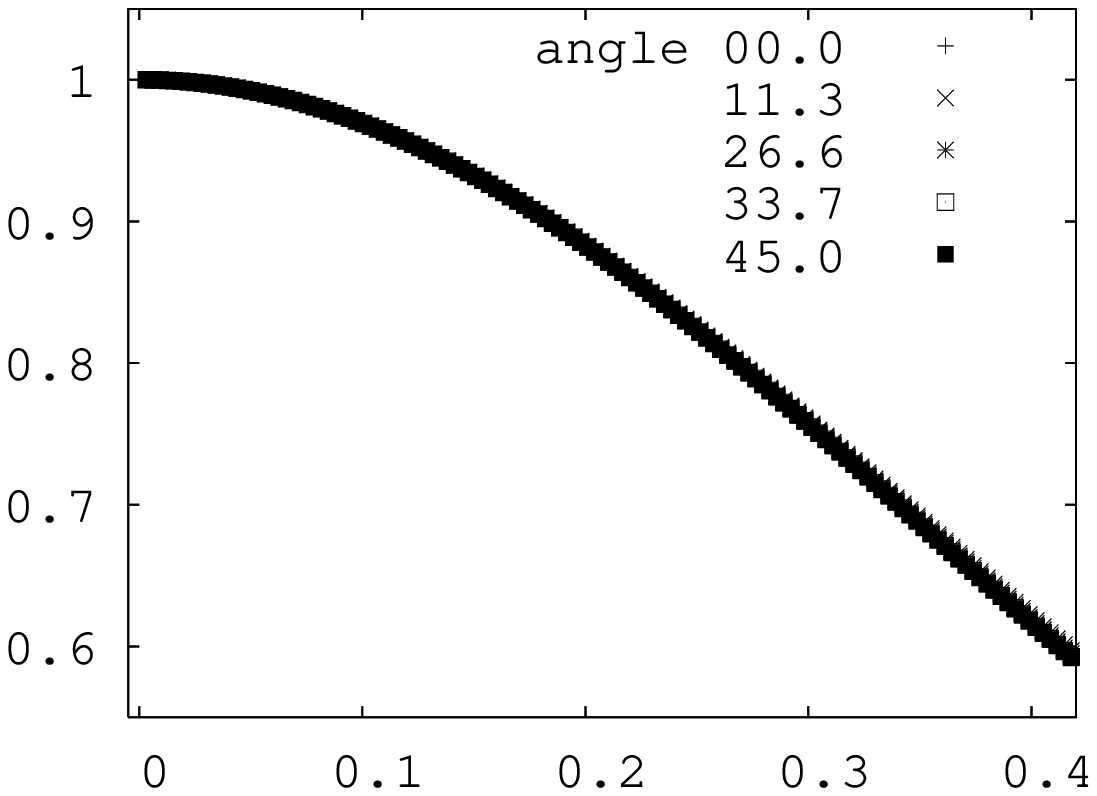}  
\includegraphics[width=.33 \textwidth] {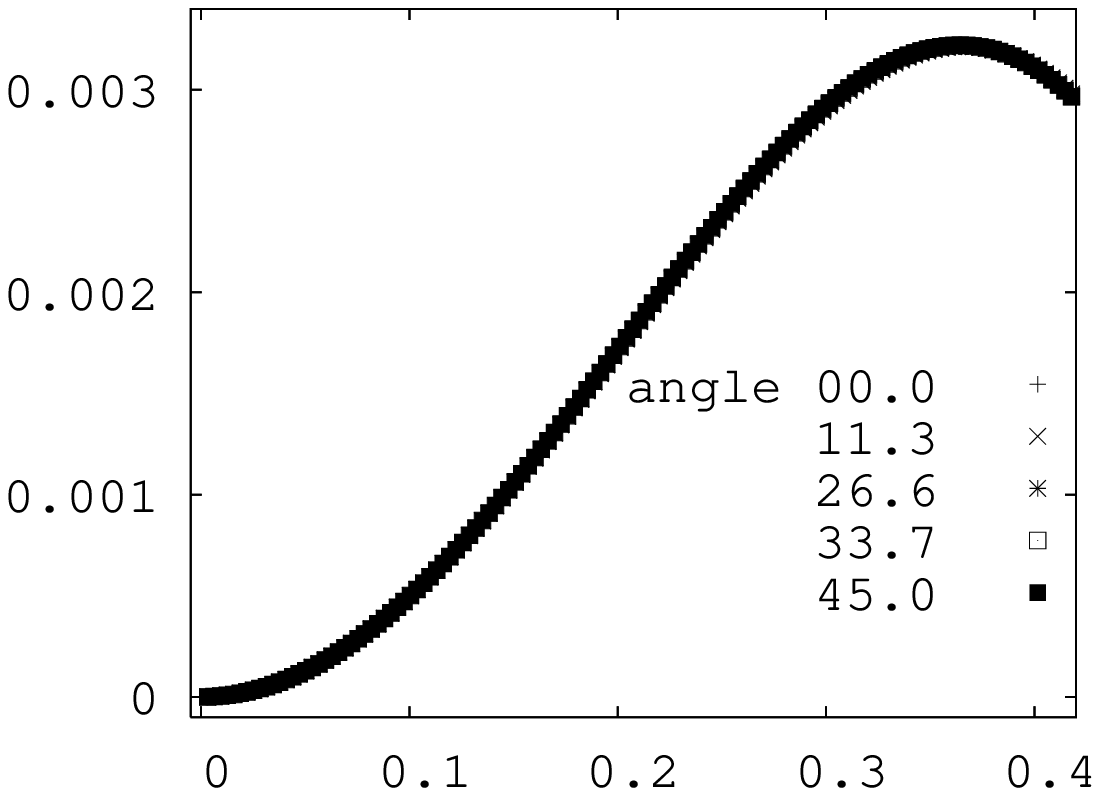} $\,$ }

\smallskip \noindent  {\bf Figure 8}. \quad ``Zero point''  experiment for the D2Q17
scheme.   
Left~:  effective viscosity  and diffusivity  $\kappa$  {\it vs} 
the wave vector  $k$ for several directions. 
Middle~:    attenuation of the sound waves.   
Right~:     $(v_{son}/c_0) - 1$       {\it vs}  $k$.  
Choice of parameters: 
$\, \alpha_2 = -619, \, $  $\, \beta_2 = -20.55, \, $      
$\, s_5 = 1.81812 ,\, $         
$\, s_{11} = 1.9230 ,\, $     
$\, s_{12} = 1.818 ,\, $      
$\, s_{17} = 1.111 .\, $      

\bigskip \monitem 
In a second case, we show that one                                                 
can reduce significantly the physical dissipations by a better tuning of the    
parameters. The associated physical parameters are given by 
\moneqstar  
Pr = 0.69817  \,, \,\, 
c_0 =  \sqrt{7\over6}  \,, \,\, 
  {{\nu}\over{\lambda \, \Delta x}}  =    0.001167       \,, \,\, 
  {{\kappa }\over{\lambda \, \Delta x}} =    0.001671    \,, \,\, 
 {{\gamma }\over{\lambda \, \Delta x}} =    0.000651    \, .      
\monendstar
The dissipation is reduced by one order of magnitude
if we refer to the previous example.                            
The results are presented in Figure~9. 
We observe that the isotropy of the waves is not rigorously 
satisfied. 
A systematic search in the space of free parameters of the model would certainly lead
to better behavior, especially in order to increase the numerical stability of the model which, as 
presented here, is not very good.                      
We refer for this approach to Xu and Sagaut \cite{XS12}.

\smallskip   \smallskip     
\centerline {
\includegraphics[width=.33 \textwidth] {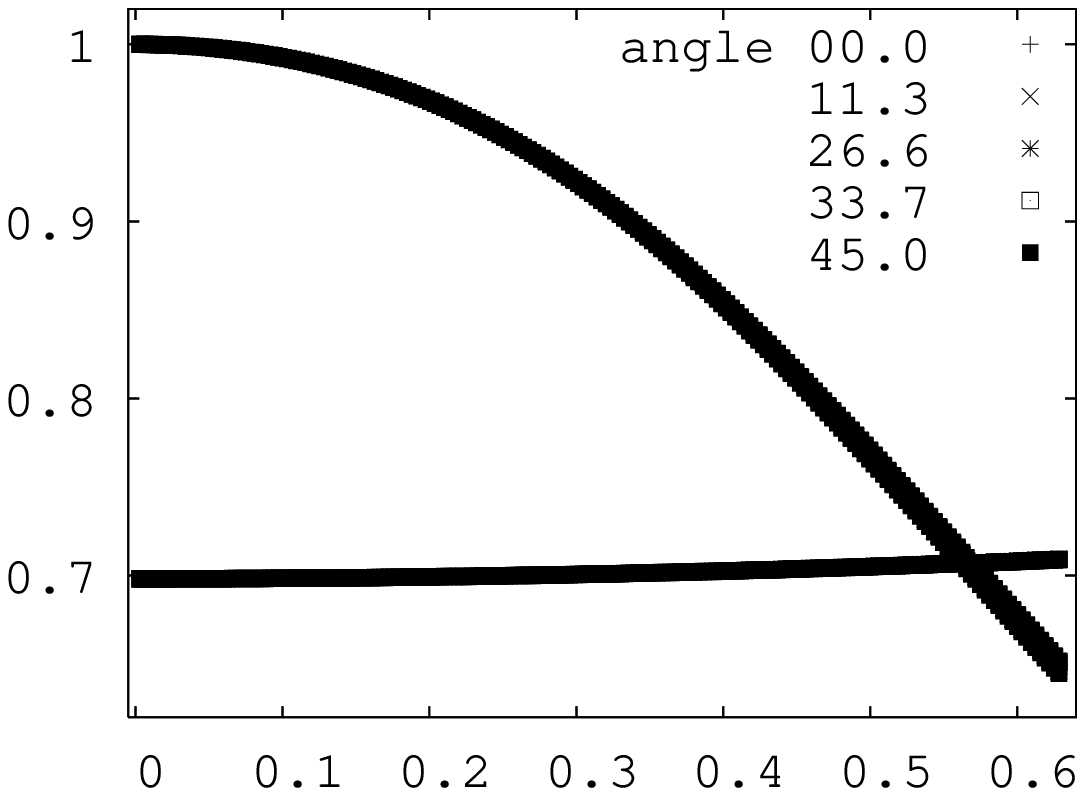}  
\includegraphics[width=.33 \textwidth] {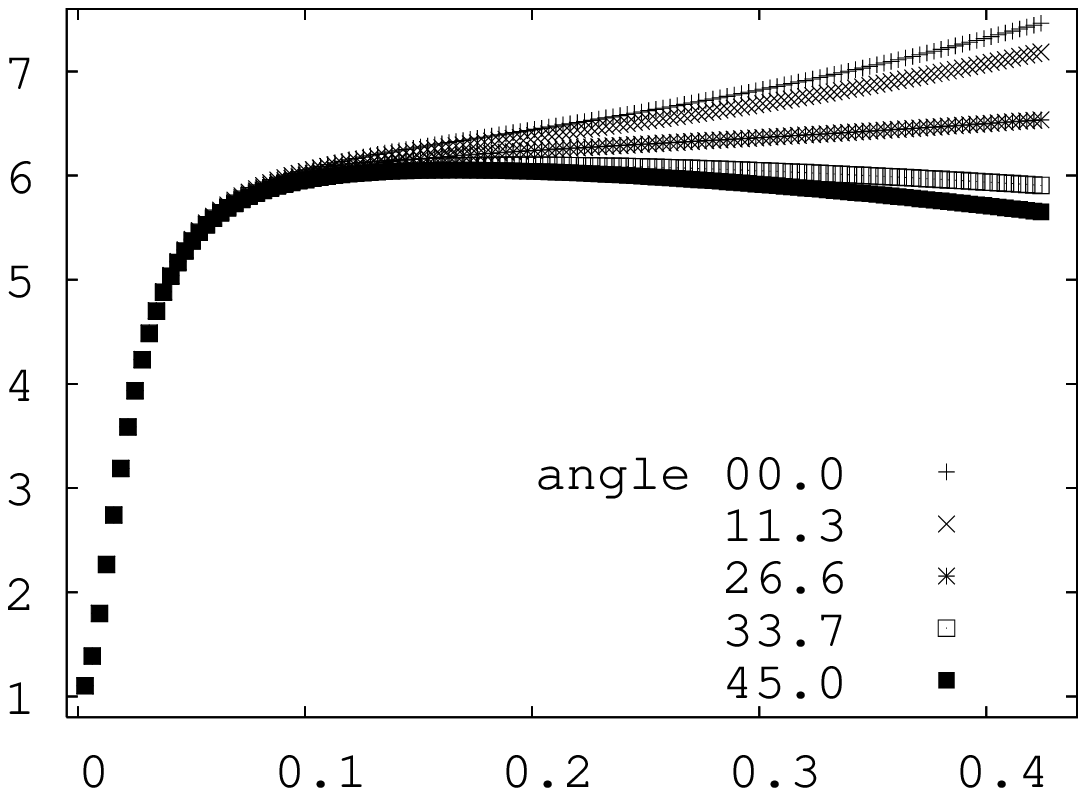}  
\includegraphics[width=.33 \textwidth] {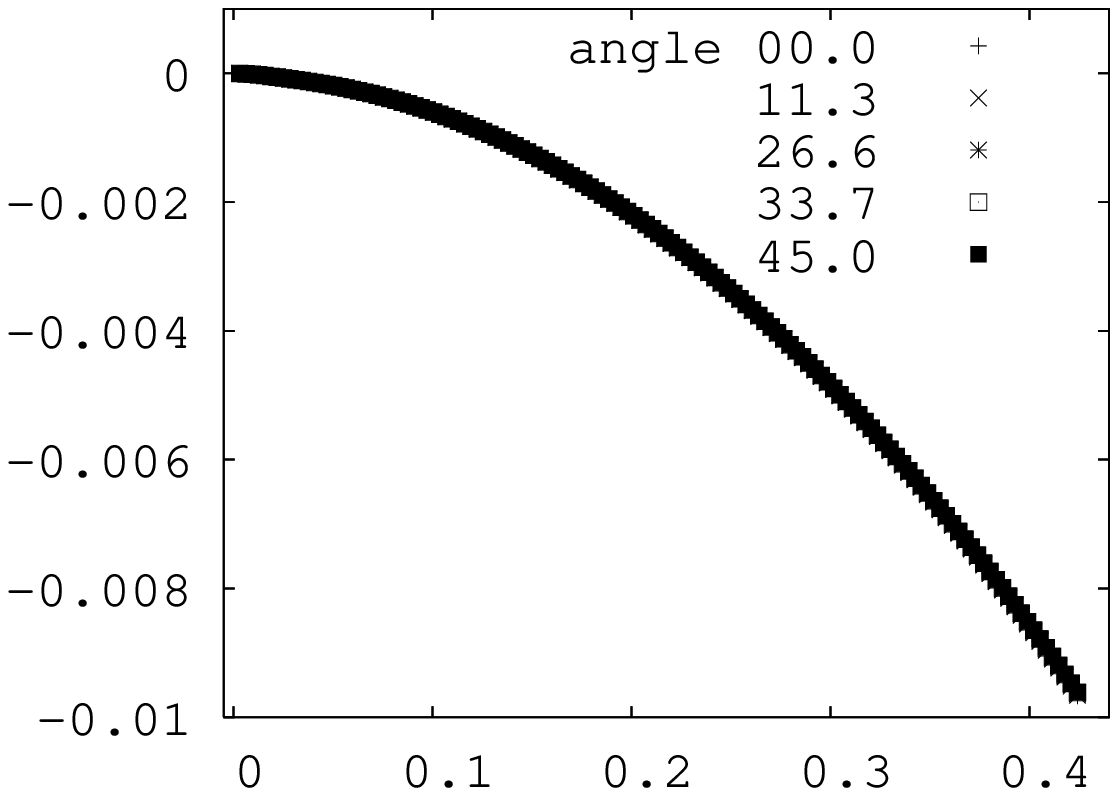}}  

\smallskip \noindent  {\bf Figure 9}. \quad ``Zero point''  experiment for the D2Q17
scheme.   
Left~:  effective viscosity  and diffusivity  $\kappa$  {\it vs} 
the wave vector  $k$ for several directions. 
Middle~:    attenuation of the sound waves.   
Right~:     $(v_{son}/c_0) - 1$       {\it vs}  $k$.   
Choice of parameters:
 $\, \alpha_2 = -641.17  \, $ 
 $\, \beta_2 = -21.01933,  \, $   
$\, s_5 = 1.9920 ,\, $         
$\, s_{11} = 1.9230 ,\, $     
$\, s_{12} = 1.818 ,\, $      
$\, s_{17} = 1.25 .\,  $      

\newpage

 \bigskip \bigskip   \noindent {\bf \large 8) \quad  Conclusion  }  

\noindent 
We have considered the 
    problem of  ``energy conserving''  lattice Boltzmann models. 
Not completely satisfying results   were proposed in the literature
with the classic  version of D2Q13 LB scheme \cite{LL03} even                
for very elementary situations as a  shear wave and  diffusive wave.               
We have added two new ideas :  add     nonlinear terms and 
remove the ``spurious coupling'' with a fourth order  analysis  
of the equivalent partial equivalent equation.
More precisely, our theoretical analysis is founded of the 
knowledge of the partial equivalent equations of the lattice 
Boltzmann scheme at several orders of precision. 
At the first order the linear nondissipative advective acoustics suggest
which nonlinear terms should be included in the equilibrium values         
 of the second order moments and the third order heat flux.          
At the second order the  linear dissipative advective acoustics 
establish general relations for the viscosity and diffusivity from the necessary isotropic     
behavior of the LBE model leads to constraint on the linear dependence of higher order moments.   
It is possible to enforce  Galilean invariance 
at first order accuracy for shear, thermal and acoustic waves.                        
The analysis of classical acoustics allows the computation of parameters 
that are compatible with isotropic waves. 
    Satisfactory results are shown for the shear wave     
for three versions of the lattice Boltzmann model considered here.
%
%
%
This breakthrough has to be confirmed  
for other test cases, lattice Boltzmann  models and higher dimensions~!

\bigskip  \bigskip   \noindent {\bf \large Acknowledgments}   

 \noindent     
The authors  thank  Nikolaos Prasianakis 
for his suggestion  to  incorporate  the D2Q9 lattice Boltzmann scheme         
in the framework of this contribution.

\bigskip \smallskip   \smallskip                     
  \centerline { \includegraphics[width=.25 \textwidth] {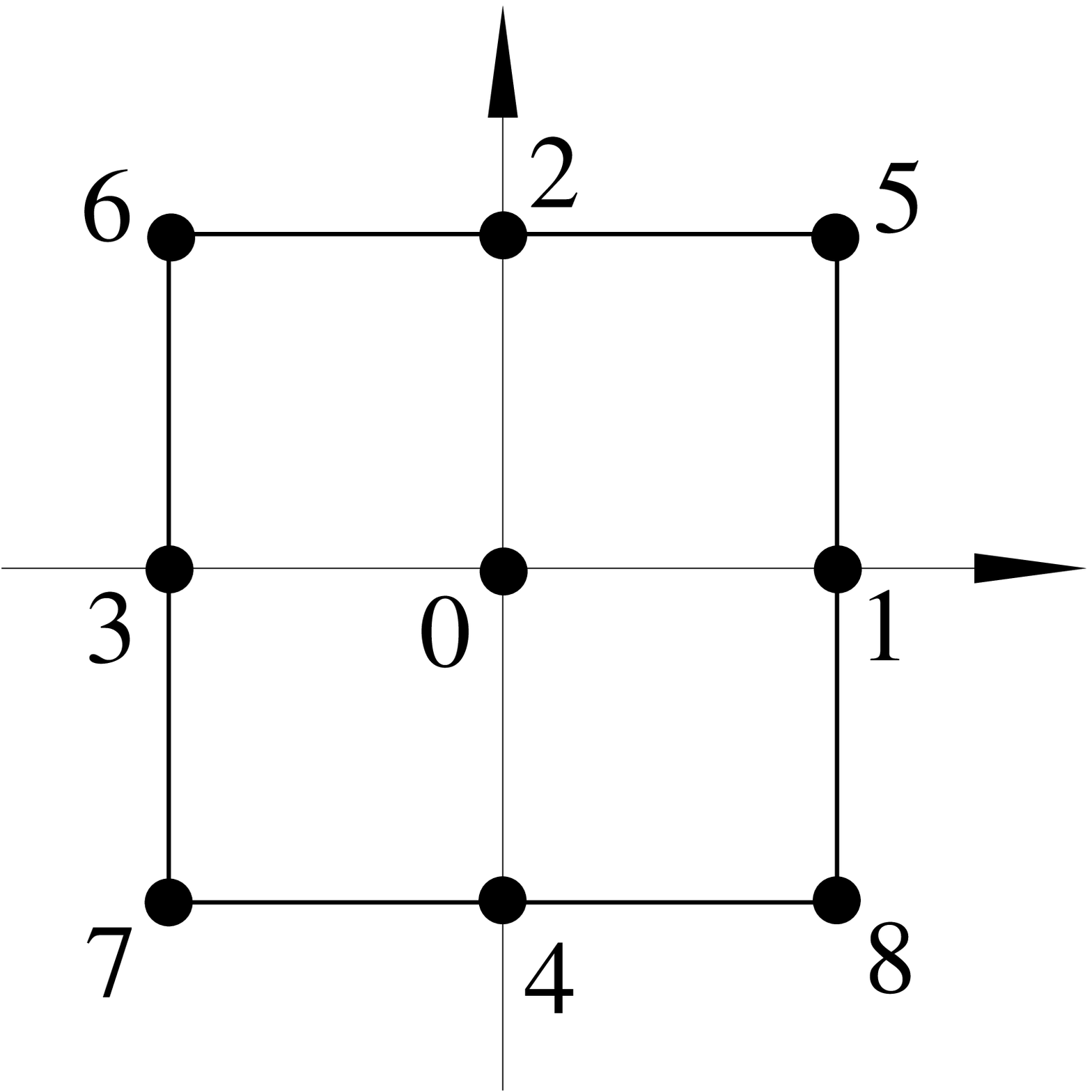}}

\smallskip \noindent  {\bf Figure 10}. \quad 
Stencil of 9  velocities for the D2Q9 lattice Boltzmann scheme.  
\smallskip \smallskip 
%

\bigskip   \noindent {\bf  \large   Annex - 1. \quad  D2Q9 lattice Boltzmann
  scheme  } 

\smallskip   \noindent 
The velocity set $ \, v_j \,$ for $ \, 0 \leq j \leq q-1  \,$ 
of a DdQq lattice Boltzmann scheme is given by the general
relation
\moneq \label{vit-schema} 
v_j \,=\, \xi_j \, \lambda \,,     
\monend 
where $\, \lambda \, $ is some   scale velocity.  
For the D2Q9  scheme \cite{LL00}  illustrated in  Figure~10, 
the $ \, \xi_j$'s of relation  (\ref{vit-schema})   are given by the expressions 
\moneq \label{vit-d2q9} \left\{ \begin{array} [c]{l}  
\displaystyle   
\xi_0 = (0,\, 0)\,, \,\, 
\xi_1 = (1,\, 0)\,, \,\,  
\xi_2 = (0,\, 1)\,, \,\, 
\xi_3 = (-1,\, 0)\,, \,\,  
\xi_4 = (0,\, -1)\,,    \\ \displaystyle 
\xi_5 = (1,\, 1)\,, \,\, 
\xi_6 = (-1,\, 1)\,, \,\,  
\xi_7 = (-1,\, -1)\,, \,\, 
\xi_8 = (1,\, -1)\,.   
\end{array} \right.     \monend

\noindent
The d'Humi\`eres  moments \cite{DdH92}   are defined with the help 
of a family $\, p_k \,$ ($0 \leq k \leq q-1$)   of  two variables polynomials. 
They are ordered by increasing degree.  
For the D2Q9 scheme, 
$ \, p_j \in   {\cal P}_{\rm D2Q9} \, $ with  
%
\moneq   \label{polynomes-d2q9} 
 \left\{   \begin{array} [l]{l}  \displaystyle  
 p_0 = 1  \,, \quad  
 p_1 = X \,, \quad 
 p_2 = Y   \,, \quad 
 p_3 =   -4 \, \lambda^2  \,+\, 3 \, (X^2+Y^2)   \\ \displaystyle 
 p_4 = X^2-Y^2   \,, \quad    
 p_5 =   X \,  Y  \\ \displaystyle 
 p_6 =  X \, \big( -5 \, \lambda^2  + 3 \, ( X^2+Y^2) \big)  \,, \,\,  
 p_7 =     Y \,  \, \big( -5 \, \lambda^2  + 3 \, ( X^2+Y^2) \big)  \\ \displaystyle 
 p_8 = 4 \,  \lambda^4 \,-\, {{21}\over{2}} \,  \lambda^2 \, (X^2+Y^2) 
\,+\,   {{9}\over{2}} \,  (X^2+Y^2)^2   \, . 
\end{array} \right.   \,.    \monend   
The coefficients of the matrix $\, M \,$ 
 are simply given by  nodal values in the velocity space: 
\moneq   \label{matrice-M} 
M_{kj} = p_k(v_j) \,, \qquad 0 \leq \, j \,, \,\, k   \, \leq  q-1 \, . 
\monend 
The moments $ \, m_k \, $ for $ \, 0 \leq k \leq q-1 \,$ are 
defined with the help of this matrix:
\moneq   \label{moments} 
m_k \,=\, \sum_j  \, M_{kj} \, f_j \,, \qquad 0 \leq \, k   \, \leq  q-1 \, . 
\monend 
The moments defined by the relations (\ref{polynomes-d2q9}) and (\ref{matrice-M})
are, due to (\ref{relax}), 
the eigenvectors of the relaxation operator of the Boltzmann 
equation with a finite number of velocities, as noticed in \cite{Du10}. 
In particular in this contribution, 
\moneq  \label{basic-moments}     \left\{ \begin{array} [c]{l}  
\displaystyle    
\rho \equiv \sum_j \, f_j \,, \,\,\, 
j_x \equiv \sum_j \, v_j^x \, f_j \,, \,\,\, 
j_y \equiv \sum_j \, v_j^y \, f_j \,, \,\,\, 
E  \equiv \sum_j \, p_3(v_j)  \, f_j \,, \\ \displaystyle 
XX  \equiv \sum_j \, p_4(v_j)  \, f_j \,, \,\,\, 
XY  \equiv \sum_j \, p_5(v_j)  \, f_j \,,  \\ \displaystyle  
q_x  \equiv \sum_j \, p_6(v_j)  \, f_j \,, \,\,\, 
q_y  \equiv \sum_j \, p_7(v_j)  \, f_j \,.   
\end{array} \right.    \monend 
We observe that due to the orthogonalization procedure, 
 the ``numerical''  total energy $\, E \,$ 
proposed at relations (\ref{basic-moments}) and effectively used 
in the simulations is related to 
the ``physical''   total energy $\, \varepsilon \,$ introduced 
in (\ref{var-conserv})  as the fourth conserved variable 
thanks to the relation
\moneq  \label{energie-E-epsil-d2q9}   
 E = 6 \,   \varepsilon - 4 \, \lambda^2 \,  \rho   \, .  
\monend 
In a first approach, we choose  the equilibria for the nonconserved moments as follows~: 
\moneq   \label{first-equil-D2Q9}   
XX^{\rm eq}  \,= \, 0  \,,\quad  XY^{\rm eq}  \,= \, 0 \,,\quad  
 {\bf    q}^{\rm eq}  \,= c_1 \,  \lambda^2 \,  {\bf j}  \,,\quad  
E_2^{\rm eq}  \,= \, \alpha_2 \,  \lambda^4 \, \rho +   \beta_2 \,  \lambda^2 \, E    \, .   \monend  
The coefficients $ \, s_k \,$ that determine the relaxation of the d'Humi\`eres moments 
\moneq   \label{relax} 
m_k^* \,=\, m_k \,+\, s_k \, \big( m_k^{\rm eq} - m_k \big) 
  \monend  
are defined  according to 
\moneq   \label{sk-d2q9}  
s_{XX} \equiv s_5 \,,\quad  
s_{XY} \equiv s_5 \,,\quad  
s_{qx} \equiv s_7 \,,\quad
s_{qy} \equiv s_7 \,,\quad   
s_{E2} \equiv s_{9} \,.
 \monend    
We set also 
\moneq    \label{sigma-s}  
\sigma_k \equiv {{1}\over{s_k}} - {1\over2} \, . 
  \monend 
Recall that the distribution $ \, f^*_j \, $ of particles after relaxation 
is defined from the moments $ \, m \, $ and the inversible matrix $M$ according to 
\moneq    \label{f-star} 
f^*_j \,=\,  \sum_k  \, \big( M^{-1} \big)_{jk} \,\, m_k^*  \, . 
  \monend 
%

\smallskip   \smallskip                     
  \centerline { \includegraphics[width=.45 \textwidth] {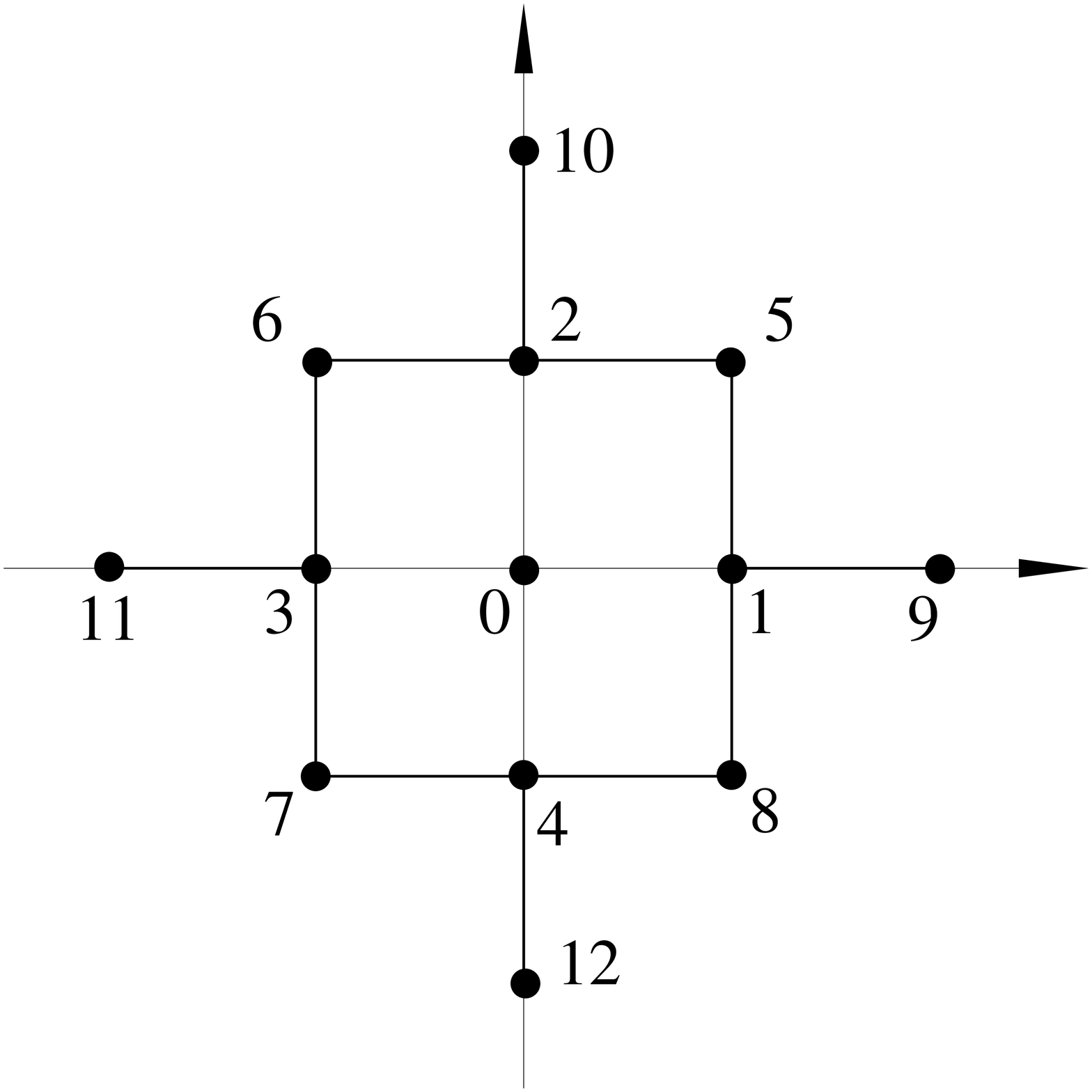}}

\smallskip \noindent  {\bf Figure 11}. \quad 
Stencil of 13 velocities for the D2Q13 lattice Boltzmann scheme.  
\smallskip \smallskip 

\bigskip \bigskip  \noindent {\bf  \large   Annex - 2. \quad  D2Q13 lattice Boltzmann
  scheme  } 

\smallskip   \noindent 
For the D2Q13 scheme \cite{Qi93, WB96}  illustrated in  Figure~11, 
the first nine $ \, \xi_j$'s of relation (\ref{vit-schema}) are the one 
given at the relation    (\ref{vit-d2q9}). The last four are  
\moneq \label{vit-d2q13}  
\xi_9 = (2,\, 0)\,, \,\,  
\xi_{10} = (0,\, 2)\,, \,\, 
\xi_{11} = (-2,\, 0)\,, \,\,  
\xi_{12} = (0,\, -2)\,, \,\,     
\monend 
The family  $\, {\cal P}_{\rm D2Q13} \,$  
of two-variable polynomials that define the moments according to 
relation (\ref{matrice-M}) are detailed as follows:  
\moneq   \label{polynomes-d2q13} 
 \left\{   \begin{array} [l]{l}  \displaystyle  
 p_0 = 1  \,, \quad  
 p_1 = X \,, \quad 
 p_2 = Y   \,, \quad 
 p_3 =   -28 \,+\, 13 \, (X^2+Y^2)   \\ \displaystyle 
 p_4 = X^2-Y^2   \,, \quad    
 p_5 =   X \,  Y  \\ \displaystyle 
 p_6 =  X \, ( -3 \, \lambda^2  + X^2+Y^2 )  \,, \,\,  
 p_7 =     Y \, (-3 \,  \lambda^2  + X^2+Y^2 )  \\ \displaystyle 
 p_8 = X \, \Big(\frac{101}{6} \lambda^4 - \frac{63}{4} \lambda^2 \, (X^2+Y^2)
  +\frac{35}{12} \ (X^2+Y^2)^2  \Big)   
 \\ \displaystyle   \vspace{-.5cm}  ~  \\ \displaystyle
 p_9 =    Y \, \Big(\frac{101}{6} \lambda^4 - \frac{63}{4} \lambda^2 \, (X^2+Y^2)  
+\frac{35}{12} \ (X^2+Y^2)^2  \Big) 
 \\ \displaystyle   \vspace{-.5cm}  ~  \\ \displaystyle
 p_{10} =      140 \, \lambda^4  - \frac{361}{2}  \lambda^2 \, (X^2+Y^2) 
+ \frac{77}{2} \, (X^2+Y^2)^2   \\ \displaystyle   \vspace{-.5cm}  ~  \\ \displaystyle
 p_{11} = -12  \,  \lambda^6 +   \frac{581}{12}   \lambda^4 \,  (X^2+Y^2) 
- \frac{273}{8} \lambda^2 \,  (X^2+Y^2)^2 + \frac{137}{24}  \,  (X^2+Y^2)^3  \\ \displaystyle
 p_{12} =   (X^2-Y^2) \ \Big(-\frac{65}{12} \lambda^2 +\frac{17}{12}\ (X^2+Y^2) \Big)  \, . 
\end{array} \right.   \,.    \monend   
The moments $ \, m_k \equiv \sum p_k(v_j) \, f_j \,$ have 
usual names given in (\ref{basic-moments}) and for the D2Q13 scheme 
by the complementary relations 
\moneq   \label{moments-d2q13} 
m_8 \,\equiv \, r_x \,, \quad   
m_9 \,\equiv \, r_y \,, \quad   
m_{10} \,\equiv \, E_2  \,, \quad   
m_{11} \,\equiv \, E_3  \,, \quad   
m_{12} \,\equiv \, XX_e  \,. \,   
\monend   
We observe also for this scheme that   the ``numerical'' total  energy $\, E \,$ 
 is a simple function of the  ``physical'' total energy $\, \varepsilon \,$.
We have  
\moneq  \label{energie-E-epsil-d2q13}   
E      = 26 \, \varepsilon -  28 \, \lambda^2 \, \rho  
\monend 
%
%
In the  first approach presented in the introduction, 
we choose the equilibria for the nonconserved moments as follows~: 
\moneq   \label{first-equil-D2Q13}  \left\{ \begin{array} [c]{l} \displaystyle 
XX^{\rm eq}  \,= \, 0  \,,\quad  XY^{\rm eq}  \,= \, 0  \,,\quad
 {\bf    q}^{\rm eq}  \,= c_1 \,  \lambda^2 \,  {\bf j}  \,,\quad   
 {\bf    r}^{\rm eq}  \,= c_2 \,  \lambda^4 \,  {\bf j}  \,,   \\ \displaystyle    
E_2^{\rm eq}  \,= \, \alpha_2 \,  \lambda^4 \, \rho +   \beta_2 \,  \lambda^2 \, E     \,,\quad 
E_3^{\rm eq}  \,= \, \alpha_3 \,  \lambda^6 \, \rho +   \beta_3 \,  \lambda^4 \,  E   \,,\quad 
XX_e^{\rm eq}  \,= \, 0   \,  . 
\end{array} \right. \monend  
The relaxation rates $ \, s_k \,$ that determine the relaxation (\ref{relax})  
of the moments are associated according to 
\moneq   \label{sk-d2q13}  \left\{ \begin{array} [c]{l} \displaystyle 
s_{XX} \equiv s_5 \,,\quad  
s_{XY} \equiv s_5 \,,\quad  
s_{qx} \equiv s_7 \,,\quad
s_{qy} \equiv s_7 \,,\quad  
s_{rx} \equiv s_9 \,,\quad
s_{ry} \equiv s_9 \,,  \\ \displaystyle  
s_{E2} \equiv s_{11} \,,\quad
s_{E3} \equiv s_{12}   \,,\quad
s_{XXe} \equiv  s_{13} \,.  
\end{array} \right. \monend    
%
%
The coefficient $\, c_1 \,$ is related to the sound velocity $ \, c_0 \,$ 
according to 
\moneq   \label{c1-d2q13}
c_1 = 2 \, c_0^2 - 3 \,   .  
\monend

\bigskip \bigskip  \noindent {\bf  \large   Annex - 3. \quad  D2Q17 lattice Boltzmann
  scheme  } 
  
\smallskip   \noindent 
For the D2Q17 scheme 
 illustrated in  Figure~12, 
the thirteen $ \, \xi_j$'s of relation (\ref{vit-schema}) are those    
given at the relation    (\ref{vit-d2q13}). The last four are  
\moneq \label{vit-d2q17}  
\xi_{13} = (2,\, 2)\,, \,\,  
\xi_{14} = (-2,\, 2)\,, \,\, 
\xi_{15} = (-2,\, -2)\,, \,\,  
\xi_{16} = (2,\, -2) \,, .      
\monend 
In a way analogous to (\ref{polynomes-d2q13}), 
the  two-variable polynomials   family  $\, {\cal P}_{\rm D2Q17} \,$   
 are given according to:  
\moneq   \label{polynomes-d2q17}  
 \left\{   \begin{array} [l]{l}  \displaystyle  
 p_0 = 1  \,, \quad  
 p_1 = X \,, \quad 
 p_2 = Y   \,, \quad  
 p_3 =  -60 \,+\, 17\, (X^2+Y^2)   \\ \displaystyle  
 p_4 = X^2-Y^2   \,, \quad    
 p_5 =   X \,  Y  \\ \displaystyle 
 p_6 =   X \, (-17 \, \lambda^2 + 3 \, (X^2+Y^2))   \,,\quad  
 p_7 =    Y \,   (-17 \, \lambda^2 + 3 \, (X^2+Y^2))  \\ \displaystyle 
\quad  {\rm p}_r =  \frac{47}{6} \lambda^4 - \frac{17}{4} \lambda^2\, (X^2+Y^2)
  +\frac{5}{12} \ (X^2+Y^2)^2  \\ \displaystyle 
 p_8 =  X \, {\rm p}_r \,,\quad   p_9 =  Y   \, {\rm p}_r   \\ \displaystyle 
\quad  {\rm p}_{\tau}  = -\frac{7429}{42} \lambda^6 +  \frac{1565}{8} \lambda^4  (X^2+Y^2) 
- \frac{2635}{48}  \lambda^2 (X^2+Y^2)^2 + \frac{465}{112}  (X^2+Y^2)^3  \\ \displaystyle 
 p_{10} = X \,  {\rm p}_{\tau}    \,,\quad    p_{11} =   Y \,   {\rm p}_{\tau}   \\ \displaystyle  
 p_{12} =  (X^2-Y^2) \ \Big( -\frac{65}{12}  \lambda^2   
+ \frac{17}{12}\ (X^2+Y^2) \Big)  
 \\ \displaystyle   \vspace{-.5cm}  ~  \\ \displaystyle
 p_{13} =    \,X \, Y \,  \Big( -\frac{65}{12}  \lambda^2 + \frac{17}{24}  (X^2+Y^2)  \Big) 
 \\ \displaystyle   \vspace{-.5cm}  ~  \\ \displaystyle   
 p_{14} =  620 \, \lambda^4  -  \frac{969}{2}  \lambda^2 \ (X^2+Y^2) 
 + \frac{109}{2} \ (X^2+Y^2)^2     
 \\ \displaystyle   \vspace{-.5cm}  ~  \\ \displaystyle   
 p_{15} =    -16740 \, \lambda^6 +  \frac{330361}{12} \lambda^4 \ (X^2+Y^2) 
- \frac{74485}{8}  \lambda^2 \    (X^2+Y^2)^2  + \cdots  \\ \displaystyle  \hfill 
+ \frac{18445}{24}  \ (X^2+Y^2)^3    
 \\ \displaystyle   \vspace{-.8cm}  ~  \\ \displaystyle  
 p_{16} =  84 \, \lambda^8 - \frac{24055}{56}  \lambda^6  (X^2+Y^2) 
+  \frac{35425}{96}  \lambda^4 \ (X^2+Y^2)^2   + \cdots  
 \\ \displaystyle   \vspace{-.5cm}  ~  \\ \displaystyle   
 \,   \qquad \qquad \qquad  \qquad \qquad \qquad  
- \frac{6035}{64}  \lambda^2 \ (X^2+Y^2)^3  
+ \frac{9193}{1344}  (X^2+Y^2)^4   \, . \,  ~ 
\end{array} \right.      \monend   
The first moments are precise  at  the relations 
 (\ref{basic-moments}). The 
new moments introduced with the D2Q17 scheme 
with the help of  relations (\ref{matrice-M}) and (\ref{polynomes-d2q17}) 
are 
\moneq  \label{moments-d2q17}  
 \left\{   \begin{array} [l]{l}  \displaystyle  
m_{8} \,\equiv \, r_x  \,, \quad   
m_{9} \,\equiv \, r_y  \,, \quad   
m_{10} \,\equiv \, \tau_x   \,, \quad   
m_{11} \,\equiv \, \tau_y    \,, \quad   
m_{12} \,\equiv \, XX_e \,,    \\ \displaystyle      
m_{13} \,\equiv \, XY_e    \,, \quad 
m_{14} \,\equiv \, E_2    \,, \quad      
m_{15} \,\equiv \, E_3    \,, \quad      
m_{16} \,\equiv \, E_4    \,.    
\end{array} \right.      \monend   
We observe between the   ``numerical'' and   ``physical'' total  energies
a relation very analogous to (\ref{energie-E-epsil-d2q9}) 
and (\ref{energie-E-epsil-d2q13}).    
We have  for the D2Q17 lattice Boltzmann scheme:
\moneq  \label{energie-E-epsil-d2q17}   
E      = 34 \, \varepsilon -  60 \, \lambda^2 \, \rho  
\monend 
Due to natural isotropy conditions, the $\, \sigma$'s coefficients defined 
by (\ref{sigma-s}), satisfy the relations 
\moneq  \label{sigma-d2q17}   
\sigma_4 = \sigma_5 \,, \quad 
\sigma_6 = \sigma_7 \,, \quad 
\sigma_8 = \sigma_9 \,, \quad 
\sigma_{10} = \sigma_{11} \,, \quad 
\sigma_{12} = \sigma_{13} \,. 
\monend 
%
 
\smallskip   \bigskip       \bigskip              
  \centerline { \includegraphics[width=.40 \textwidth] {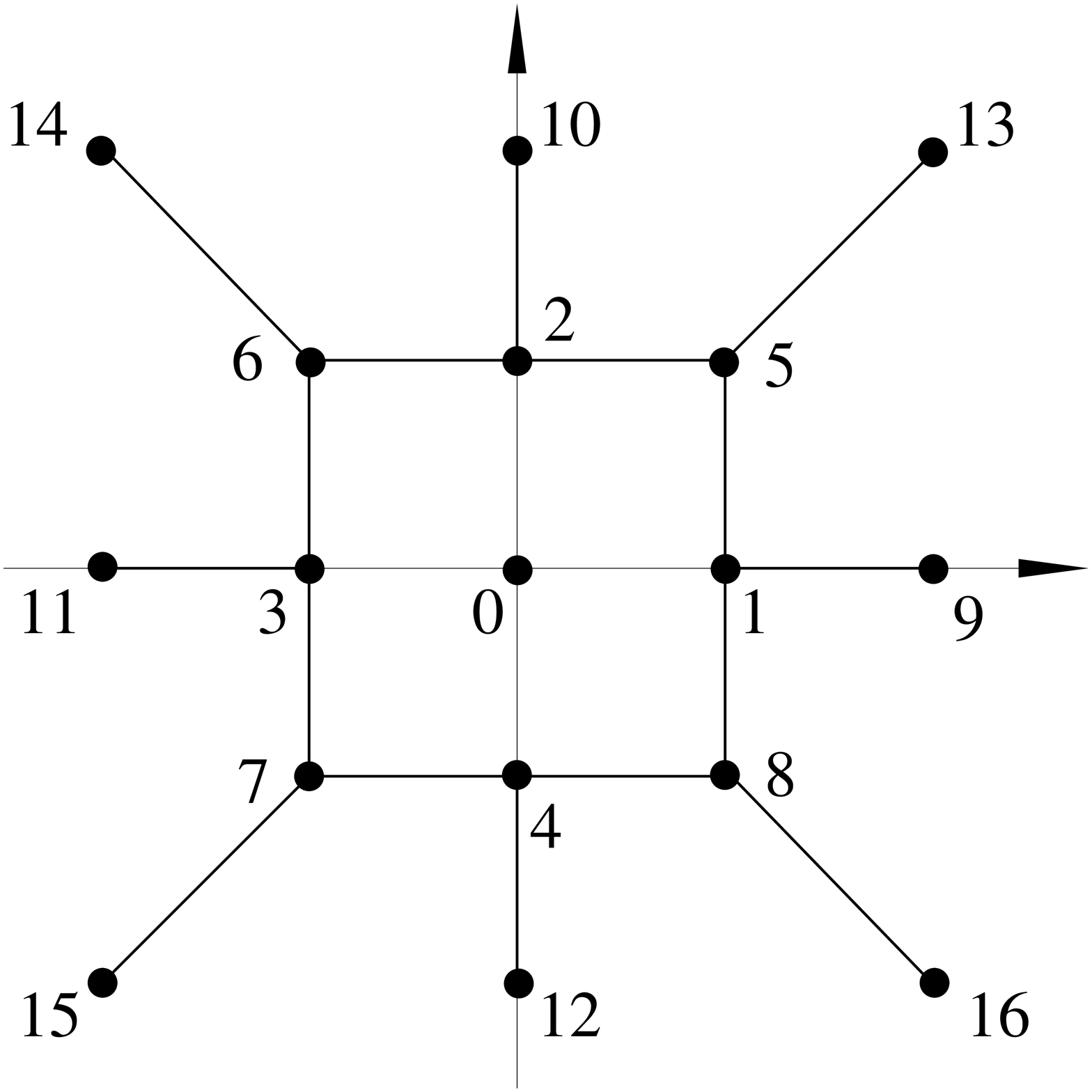}} 

\smallskip \noindent  {\bf Figure 12}. \quad 
Stencil of 17 velocities for the D2Q17 lattice Boltzmann scheme.   
\smallskip \smallskip 

\bigskip \bigskip  \bigskip 

  \end{document}